\documentclass[AMA,STIX1COL]{WileyNJD-v2}
\usepackage{moreverb}

\usepackage{amsmath,amsthm,amsfonts,mathtools,bm}
\usepackage{dsfont}
\usepackage{tikz,pgfplots}
\usetikzlibrary{calc,patterns,
                 decorations.pathmorphing,
                 decorations.markings}
\usepackage{subcaption}
\pgfplotsset{compat=newest}
\usepackage{float,indentfirst,algorithm}
\usepackage{pgfplots}
\usepackage{booktabs}
\usepackage{hyperref}
\usepackage{bbm}
\usepackage{makecell}

\usepackage[capitalise]{cleveref}

\newcommand{\mean}{m}
\newcommand{\pmean}{\widehat{m}}
\newcommand{\pp}{\widehat{\theta}}
\newcommand{\py}{\widehat{y}}
\newcommand{\ppy}{\widehat{y}}
\newcommand{\Cov}{C}
\newcommand{\pCov}{\widehat{C}}

\newcommand{\N}{\mathcal{N}}
\newcommand{\G}{\mathcal{G}}

\newcommand{\E}{\mathbb{E}}
\newcommand{\I}{\mathds{I}}
\newcommand{\R}{\mathbb{R}}

\newcommand{\Cv}{\mathcal{C}}

\DeclareMathOperator*{\argmin}{arg\,min}

\newcommand\BibTeX{{\rmfamily B\kern-.05em \textsc{i\kern-.025em b}\kern-.08em
T\kern-.1667em\lower.7ex\hbox{E}\kern-.125emX}}

\articletype{Article Type}%

\received{<day> <Month>, <year>}
\revised{<day> <Month>, <year>}
\accepted{<day> <Month>, <year>}

\begin{document}

\title{Bayesian Calibration for Large-Scale Fluid Structure Interaction Problems Under Embedded/Immersed Boundary Framework}

\author[1]{Shunxiang Cao}

\author[2]{Daniel Zhengyu Huang*}

\authormark{Cao and Huang}

\address[1]{\orgdiv{Department of Mechanical and Civil Engineering}, \orgname{California Institute of Technology}, \orgaddress{\state{CA}, \country{USA}}}

\address[2]{\orgdiv{Department of Environmental Science and Engineering, Department of Computing Mathematical Sciences}, \orgname{California Institute of Technology}, \orgaddress{\state{CA}, \country{USA}}}

\corres{*Daniel Zhengyu Huang, Department of Environmental Science and Engineering, Department of Computing Mathematical Sciences, California Institute of Technology.\\
 \email{dzhuang@caltech.edu}}

% \presentaddress{}

\abstract[Abstract]{

Bayesian calibration is widely used for inverse analysis and uncertainty analysis for complex systems in the presence of both computer models and observation data.
In the present work, we focus on large-scale fluid-structure interaction systems characterized by large structural deformations.
Numerical methods to solve these problems, including embedded/immersed boundary methods, are typically not differentiable and lack smoothness.
We propose a framework that is built on unscented Kalman filter/inversion to efficiently calibrate and provide uncertainty estimations of such complicated models with noisy observation data.
The approach is derivative-free and non-intrusive, and is of particular value for the forward model that is computationally expensive and provided as a black box which is impractical to differentiate.
The framework is demonstrated and validated by successfully calibrating the model parameters of a piston problem and identifying the damage field of an aircraft wing under transonic buffeting.
}

\keywords{Fluid-structure interaction; Derivative-free optimization; Bayesian calibration; Embedded/Immersed boundary method; Uncertainty quantification}

\maketitle

\footnotetext{\textbf{Abbreviations:} FSI, fluid structure interaction; EBM, embedded boundary method;  UKI, unscented Kalman inversion; CFD, computational fluid dynamics; AMR, adaptive mesh refinement; MCMC Markov chain Monte Carlo; ETKI, ensemble transform Kalman inversion.}

\section{Introduction}
Fluid‐Structure Interaction (FSI) problems arise in many scientific and engineering applications including, to name only a few, aircraft aeroelasticity~\cite{geuzaine2003aeroelastic,kamakoti2004fluid,dowell2010some}, parachute inflation dynamics~\cite{huang2018simulation,huang2020modeling,avery2021computationally}, hemodynamics~\cite{hsu2014fluid,updegrove2017simvascular,bertaglia2021uncertainty}, and lithotripsy~\cite{cao2019shock,cao2021shock}.
Besides the development of mathematical models, seamless integration of observation data with these models starts to play a significant role to improve the prediction and quantify uncertainty for FSI, for example, calibration of hemodynamic model parameters to match patient data~\cite{ismail2013adjoint,arthurs2020flexible,bertaglia2021uncertainty} and structural damage detection using sensor data and a digital twin~\cite{sohn1997bayesian,tuegel2011reengineering,sotoudehnia2019new}.
The integration can be formulated as a data-based calibration problem. And major associated challenges include:
\begin{itemize}
    \item Observation data are noisy;
    \item The calibration problems might not be done for the single fluid or structure subsystem without FSI coupling, due to strong fluid-structure coupling and complex mechanism (e.g., viscoelasticity~\cite{bertaglia2021uncertainty}), which requires FSI solvers with advanced numerical methods particularly designed for the applications of interest;
    \item The FSI solvers might be given as a black box (difficult to calculate the derivative in practice), or the solvers are not differentiable due to the numerical methods (e.g., embedded boundary method and adaptive mesh refinement) or the physical nature (fracture);
    \item Each forward FSI evaluation is expensive for real-world applications.  
\end{itemize}
Therefore, an efficient non-intrusive algorithm for calibration and uncertainty quantification is highly desirable. 

The present work focuses on data-based calibration for FSI problems where the structure undergoes large displacements, large rotations and/or deformations, as well as topological changes.
%FSI applications of particular interest in the present work are problems in which the structure undergoes large displacements, large rotations and/or deformations, as well as topological changes.
A popular class of methods for solving this kind of FSI problems are the Embedded~\cite{lohner2004adaptive,wang2011algorithms,farhat2012fiver,lakshminarayan2014embedded,main2017enhanced} or Immersed Boundary Methods~\cite{peskin1972flow,fadlun2000combined,kim2001immersed,uhlmann2005immersed,choi2007immersed,taira2007immersed,tian2014fluid} (EBMs or IBMs). These methods compute fluid flow on non-body-fitted computational fluid dynamics~(CFD) meshes in which discrete representations of wet surfaces of obstacles are embedded or immersed.
There are variants under other names, including cut cell methods~\cite{ingram2003developments, schott2019monolithic,berger2012progress}, fictitious domain methods~\cite{johansen1998cartesian,uddin2014cartesian,harris2012adaptive,almgren1997cartesian,balaras2004modeling}, ghost fluid-structure methods~\cite{tseng2003ghost}, and immersed boundary–Lattice Boltzmann methods~\cite{feng2004immersed,tian2011efficient}.
However, as for data-based calibration, EBMs and IBMs might not be favorable, because they are generally not differentiable and the quantities of interest, like surface stresses and forces, require special treatments to retain smoothness~\cite{uhlmann2005immersed,yang2009smoothing, goza2016accurate,ho2020discrete,ho2021aerodynamic}.
The non-differentiation and lack of smoothness are rooted in the enforcement of fluid-structure interface conditions on a non-interface-conforming mesh.
%The non-differentiation and lack of smoothness are rooted in the enforcement of the transmission conditions on the non-conforming mesh on the fluid-structure interface.
More specifically, the stencils, which are used to evaluate discrete delta function or reconstruct fluid states at the "sharp interface", keep changing along with the moving interface. And the status of a node may switch between real fluid node and ghost fluid node as the structure moves through the fixed mesh, producing severe oscillations in the solution.
%Because the stencils, which are used to evaluate discrete delta function or reconstruct fluid states at the "sharp interface", keep changing along with the moving interface. And the status of the fluid nodes, which can be classified as real or active if they are located inside the region of the computational fluid domain occupied by the fluid and ghost or inactive otherwise, keeps switching. 
These discrete events render adjoint-based optimization approaches almost impractical for the calibration.
Moreover, EBMs and IBMs are generally combined with adaptive mesh refinement (AMR)~\cite{berger1989local,griffith2007adaptive,borker2019mesh} for better resolution of the fluid-structure interface. 
AMR "discretely" adds and removes fluid nodes for refinement and coarsening, which further complicates the differentiation.
Therefore, in the present work, we focus on derivative-free approaches for calibration.
It is worth mentioning although the calibration approach is demonstrated with the EBM, it is equally applicable to its body-fitted counterpart---arbitrary Lagrangian-Eulerian~(ALE) approach~\cite{hirt1974arbitrary,tezduyar1992new,felippa2001partitioned,kuttler2008fixed,huang2019high,huang2020high, huang2019highaiaa}, which relies on mesh motion, deformation schemes, and local remeshing~\cite{dukowicz1987accurate,tezduyar1992new,long2013fluid} to maintain mesh conformity at the fluid-structure interface.

Derivative-free Bayesian calibration or inversion~\cite{kaipio2006statistical,dashti2013bayesian} generally starts with the observation error model,
\begin{equation}
\label{eq:obs-err-model}
y  = \G(\theta) + {\eta},    
\end{equation}
where the forward operator
$\G: \R^{N_{\theta}} \to \R^{N_y}$ maps the unknown model parameter vector $\theta \in \R^{N_{\theta}}$
to the observation vector $y \in \R^{N_y}$. Specifically, the operator $\G$ represents the FSI solver with proper initial and boundary conditions.  
For a given observation, the observational noise $\eta$ is unknown, but it is assumed to statistically follow a known distribution. To be concrete we will assume
that it is drawn from a Gaussian with distribution $\N(0,\Sigma_{\eta})$.
Given a guess of the distribution of $\theta$, represented by a prior density function $p_{0}(\theta)$, Bayesian calibration aims to estimate the posterior distribution of $\theta$ that satisfies
\begin{equation}
\label{eq:post}
    p(\theta) \propto \exp\bigl(-\Phi(\theta)\bigr)p_0(\theta) \qquad \textrm{with} \quad \Phi(\theta) = \frac{1}{2}\lVert\Sigma_{\eta}^{-\frac{1}{2}}(y - \G(\theta)) \rVert^2,
\end{equation}
where $\Phi(\theta)$ denotes the misfit between the modeled and observed data. The resulting $p(\theta)$ is supposed to provide a good confidence interval of the truth parameter.
%The posterior density function $p$ of the model parameter $\theta$ is
%\begin{equation}
%\label{eq:post}
%    p(\theta) \propto \exp\bigl(-\Phi(\theta)\bigr)p_0(\theta) \qquad \textrm{with} \quad \Phi(\theta) = \frac{1}{2}\lVert\Sigma_{\eta}^{-\frac{1}{2}}(y - \G(\theta)) \rVert^2,
%\end{equation}
%where $p_0(\theta)$ represents the prior density function.
%The objective of Bayesian calibration is to estimate the posterior distribution of $\theta$, which is supposed to provide a good confidence interval of the truth parameter. 
In the present work, we assume there are enough observation data, and use improper uniform prior $p_0 \propto 1$ to let the data speak for itself. 
It is worth noticing that the maximum a posteriori (MAP) estimation for $\theta$ with the improper uniform prior corresponds to the minimizer of $\Phi(\theta)$. Mostly, the calibration is simplified and reformulated as a nonlinear least-square optimization problem to find optimal $\theta$ that satisfies
\begin{equation}
\label{eq:LSQ}
    \theta = \argmin_{\theta} \frac{1}{2}\lVert\Sigma_{\eta}^{-\frac{1}{2}}(y - \G(\theta)) \rVert^2
\end{equation}
instead of its distribution. 
In this case, the major challenge of the optimization also lies in the aforementioned fact that $\G$ is not differentiable.

Traditional methods for derivative-free Bayesian calibration to estimate the posterior distribution, such as Markov chain Monte Carlo~\cite{geyer1992practical,gelman1997weak,goodman2010ensemble,larson2019bayesian}
(MCMC), typically require many iterations—often more than $10^4$---to reach statistical convergence. 
Given that each forward run can be expensive, conducting so many runs is computationally unaffordable, rendering MCMC impractical for real-world FSI calibrations.
%Conducting so many computationally expensive forward runs is unaffordable, rendering MCMC impractical for real-world FSI calibrations. 

In the present work, we employ the Kalman inversion for the Bayesian calibration, which approximates the posterior distribution as Gaussian distribution and generally requires $O(10)$ iterations.
Kalman inversion derived from Kalman filter~\cite{kalman1960new} and its variants, including but not limited to extended Kalman filter~\cite{sorenson1985kalman}, ensemble Kalman filter~\cite{evensen1994sequential}, unscented Kalman filter~\cite{julier1995new,wan2000unscented}, and cubature Kalman filter~\cite{arasaratnam2009cubature}, which are developed to sequentially update the probability distribution of states in partially observed  dynamics.
Kalman filtering is a two-step procedure: (1) the prediction step, where the state is computed forward in time;  
(2) the analysis step, where the state and its uncertainty are corrected to take into account the observation data.
In the analysis step, Kalman filters use Gaussian ansatz to formulate Kalman gain to assimilate the observation and update the distribution.
Numerous applications of Kalman filters, including weather forecasts 
~\cite{evensen1994sequential,anderson2001ensemble,bishop2001adaptive} and guidance, navigation, and control of vehicles~\cite{sorenson1985kalman,julier1995new, arasaratnam2009cubature} demonstrate empirically that Kalman filters can not only calibrate the model predication, but also provides uncertainty information for state estimation problems.
Kalman inversion which applies Kalman filters for parameter calibration originates from~\cite{wan1997neural,wan2000unscented,gu2006ensemble,oliver2008inverse,chen2012ensemble,iglesias2013ensemble,UKI1}.
Specifically, the parameter-to-data map~(\cref{eq:obs-err-model}) is first paired with a stationary stochastic dynamical system for the parameter, and then techniques from Kalman filtering are iteratively applied to estimate the parameter given the observation data. 
Consider the following stochastic dynamical system~\cite{UKI1},
\begin{subequations}
\label{eq:dynamics}
  \begin{align}
  &\textrm{evolution:}    &&\theta_{n+1} = \theta_{n} +  \omega_{n+1}, &&\omega_{n+1} \sim \N(0,\Sigma_{\omega}) \label{eq:evolve}\\
  &\textrm{observation:}  &&y_{n+1} = \G(\theta_{n+1}) + \nu_{n+1}, &&\nu_{n+1} \sim \N(0,\Sigma_{\nu})
\end{align}
\end{subequations}
where $\theta_{n+1}$ is the unknown parameter vector at the artificial time $n+1$, and $y_{n+1} =y$ is the observation, the artificial evolution error $\omega_{n+1}$ and artificial observation error $\nu_{n+1}$ are mutually independent, zero-mean Gaussian sequences with covariances $\Sigma_{\omega}$ and $\Sigma_{\nu}$, respectively. 
It is worth distinguishing the artificial time from the real time for time-dependent problems. In the present work, each artificial time step represents a Kalman filtering iteration in which the observation $y$ may consist of the time-series observation data from a full forward simulation.
Given its non-intrusive nature, Kalman inversion has been widely used as an optimization method for parameter estimation~~\cite{iglesias2016regularizing, schillings2017analysis, schneider2017earth, schillings2018convergence, iglesias2020adaptive,chada2020iterative,gao2021bi,zhang2021assimilation}, especially for problems where the forward model is expensive and provided as a black box that is impractical to differentiate.
%Kalman inversion has been widely used as a non-intrusive optimization method for parameter estimation~~\cite{iglesias2016regularizing, schillings2017analysis, schneider2017earth, schillings2018convergence, iglesias2020adaptive,chada2020iterative,gao2021bi,zhang2021assimilation}. And it is of particular value for parameter estimation problems where the forward model is expensive and provided as a black box which is impractical to differentiate.

Specifically, unscented Kalman inversion~\cite{wan2000unscented,UKI1} (UKI) is employed in the present work, where unscented Kalman filter~\cite{julier1995new,julier1997new,wan2000unscented} is iteratively applied to the stationary stochastic dynamics (Eq.~\eqref{eq:dynamics}).
From the optimization point of view, UKI converges exponentially fast for linear calibration problems~\cite{UKI1}; UKI performs like a generalized Levenberg-Marquardt Algorithm with a smoothed data-misfit for nonlinear calibration problems. It has been demonstrated effective for handling noisy observation data and solving even chaotic problems. 
From the Bayesian point of view, 
UKI approximates the parameter distribution by a Gaussian and allows to approximate the mean and covariance of the posterior distribution. 
In the present work, we focus on using UKI for Bayesian calibration and make the following contributions:
\begin{itemize}
    \item We further develop unscented Kalman inversion methodology and establish the linear analysis about applying UKI for Bayesian calibration. 
    \item We demonstrate that UKI is able to efficiently calibrate large-scale FSI problems under embedded/immersed boundary framework.
    \end{itemize}

The remainder of this paper is organized as follows. 
In \Cref{sec:Bayesian}, Bayesian calibration and unscented Kalman inversion are introduced. 
\Cref{sec:FSI} provides an overview of the FSI under embedded boundary framework.
Numerical examples that demonstrate and validate the proposed approach are presented in~\Cref{sec:app}, involving a 1D piston problem and a challenging 3D wing damage detection problem with transonic buffeting flows.
Finally, conclusions are offered in \Cref{sec:Conclusion}.

\section{Bayesian Calibration}
\label{sec:Bayesian}
Recall that the Bayesian calibration approach starts to pair the parameter-to-data relationship encoded in \cref{eq:obs-err-model}
with the stochastic dynamical system for the parameter~(\cref{eq:dynamics}). 
A useful way to think of the procedure is through the analogy to pseudo-time stepping for solving steady-state problem.
To approximate the posterior distribution, we employ techniques
from unscented Kalman filtering to iteratively update the conditional probability density function 
$p_n$ of $\theta_n|Y_n$, where $Y_{n} := \{y_1, y_2,\cdots , y_{n}\}$ is
the observation set at artificial time $n$, until its convergence.
The probability density function $p_n$ is approximated to be Gaussian\footnote{A justification of using Gaussian to approximate the posterior distribution~\eqref{eq:post} is from the Bernstein-von Mises theorem~\cite{le2012asymptotics,van2000asymptotic,freedman1999wald,lu2017gaussian}, which states the posterior distribution becomes asymptotically a multivariate normal distribution with the increasing of data under certain regularity conditions.} $\N(m_n, C_n)$ with mean $m_n$ and covariance $C_n$ for the sake of efficiency.  
Hence, in Subsection \ref{ssec:gau} we first introduce
a Gaussian approximation algorithm, which
maps the space of Gaussian measures into itself at each step of the iteration; 
Subsection \ref{ssec:exki} shows how this algorithm can be made practical by applying a quadrature rule (i.e., unscented transform) to
evaluate certain integrals appearing in the conceptual Gaussian approximation, which leads to the UKI algorithm. In Subsection \ref{ssec:analysis}, we present theoretical and numerical analysis of the UKI in the linear setting.
\subsection{Gaussian Approximation Algorithm}
\label{ssec:gau}
At each iteration $n$, $p_n$ is updated through the prediction and
analysis steps~\cite{reich2015probabilistic,law2015data}: $p_n 
\mapsto \hat{p}_{n+1}$, and then 
$\hat{p}_{n+1} \mapsto p_{n+1}$, where $\hat{p}_{n+1}$ denotes the conditional probability density function
of $\theta_{n+1}|Y_n$. 

In the prediction step,  $\hat{p}_{n+1} = \N(\pmean_{n+1}, \pCov_{n+1})$ is also Gaussian under~\cref{eq:evolve} and satisfies
\begin{equation}
\label{eq:KF_pred}
\begin{split}
&\pmean_{n+1} = \E[\theta_{n+1}|Y_n] =  \mean_n, \qquad \pCov_{n+1} = \mathrm{Cov}[\theta_{n+1}|Y_n] = \Cov_{n} + \Sigma_{\omega}.
\end{split}
\end{equation}

In the analysis step, the joint distribution of  $\{\theta_{n+1}, y_{n+1}\}|Y_{n}$ can be approximated by a Gaussian distribution
\begin{equation}
\label{eq:KF_joint}
     \N\Bigl(
    \begin{bmatrix}
    \pmean_{n+1}\\
    \widehat{y}_{n+1}
    \end{bmatrix}, 
    \begin{bmatrix}
   \pCov_{n+1} & \pCov_{n+1}^{\theta p}\\
    {{\pCov_{n+1}}^{\theta p}}{}^{T} & \pCov_{n+1}^{pp}
    \end{bmatrix}
    \Bigr),
\end{equation}
with
\begin{equation}
\label{eq:KF_joint2}
\begin{split}
    \widehat{y}_{n+1}   =     & \E[\G(\theta_{n+1})|Y_n], \\
     \pCov_{n+1}^{\theta p} =     &  \mathrm{Cov}[\theta_{n+1}, \G(\theta_{n+1})|Y_n],\\
    \pCov_{n+1}^{p p} = &  \mathrm{Cov}[\G(\theta_{n+1})|Y_n] + \Sigma_{\nu}.
\end{split}
\end{equation}
Conditioning the Gaussian in \eqref{eq:KF_joint} to find $\theta_{n+1}|\{Y_n,y_{n+1}\}=\theta_{n+1}|Y_{n+1}$ gives the following
expressions for the mean $\mean_{n+1}$ and covariance $\Cov_{n+1}$ of the
approximation to $p_{n+1}:$
\begin{equation}
\label{eq:KF_analysis}
    \begin{split}
        \mean_{n+1} &= \pmean_{n+1} + \pCov_{n+1}^{\theta p} (\pCov_{n+1}^{p p})^{-1} (y_{n+1} - \widehat{y}_{n+1}),\\
         \Cov_{n+1} &= \pCov_{n+1} - \pCov_{n+1}^{\theta p}(\pCov_{n+1}^{p p})^{-1} {\pCov_{n+1}^{\theta p}}{}^{T}.
    \end{split}
\end{equation}%
\Cref{eq:KF_pred,eq:KF_joint,eq:KF_joint2,eq:KF_analysis} establish a conceptual algorithm for application of Gaussian approximation for Bayesian calibration.
And the integrals appearing in~\cref{eq:KF_joint2} are approximated by the unscented approach, which is detailed in the following subsection.

\subsection{Unscented Kalman Inversion}
\label{ssec:exki}
The unscented Kalman inversion uses the following unscented transform~\cite{UKI1} to evaluate \cref{eq:KF_joint2}:

\begin{definition}
Let denote Gaussian random variable $\theta \sim \N(\mean, \Cov) \in \R^{N_{\theta}}$, $2N_{\theta}+1$ symmetric sigma points are chosen deterministically:
\begin{equation}
\begin{split}
    \theta^0 = \mean \qquad \theta^j = \mean + c_j [\sqrt{\Cov}]_j \qquad \theta^{j+N_\theta} = \mean - c_j [\sqrt{\Cov}]_j\qquad (1\leq j\leq N_\theta),
\end{split}
\end{equation}
where $[\sqrt{\Cov}]_j$ is the $j$th column of the Cholesky factor of $\Cov$. The quadrature rule approximates the mean and covariance of the transformed variable $\G_i(\theta)$ as follows,  
\begin{equation}
\label{eq:ukf}
    \E[\G_i(\theta)] \approx \G_i(\theta^0)\qquad 
    \mathrm{Cov}[\G_1(\theta),\G_2(\theta)]  \approx \sum_{j=1}^{2N_{\theta}} W_j^{c} (\G_1(\theta^j) - \E\G_1(\theta))(\G_2(\theta^j) - \E\G_2(\theta))^T. 
\end{equation}
Here these constant weights are 
\begin{align*}
    &c_j = a\sqrt{N_\theta}~(j=1,\cdots,N_{\theta})
    \qquad W_j^{c} =\frac{1}{2a^2 N_\theta}~(j=1,\cdots,2N_{\theta})\qquad
     a=\min\{\sqrt{\frac{4}{N_\theta + \kappa}},  1\}.
\end{align*}
\end{definition}
The aforementioned unscented transform is different from the original unscented transform~\cite{julier2000new}.
The modification we employ here replaces the original 2nd-order approximation of the $\E[\G_i(\theta)]$ with its 1st-order counterpart, in order to
avoid negative weights.

%%%%%%%%%%%%%%%%%%%%%%%

When the unscented transform is applied to make such approximation, we obtain the unscented Kalman inversion algorithm.
The hyperparameters are chosen as 
\begin{equation}
\label{eq:hyperparameters}
    \Sigma_{\nu} = 2\Sigma_{\eta} \quad \textrm{ and } \quad \Sigma_{\omega} =  \Cov_n,
\end{equation} 
at the $n$-th iteration. This guarantees that the converged mean and covariance well approximate those of the posterior distribution~(\cref{eq:post}) with an uninformative prior for well-defined inverse problems under certain conditions (See~Subsection~\ref{ssec:analysis}). The algorithm for unscented Kalman inversion is summarized in~\cref{alg:UKI}. It is worth mentioning that in each iteration, UKI requires solving $2N_{\theta}+1$ forward evaluations parallelly, where $N_{\theta}$ is the number of unknown model parameters.

\begin{algorithm}[ht]
 \caption{Unscented Kalman Inversion}
 \label{alg:UKI}
 \begin{algorithmic}[1]
 \Function{\texttt{UKI}}{$\G,\,\mean_0,\, \Cov_0,\, y,\, \Sigma_{\eta}$}
 \While{$n \leq N_{\textrm{max}}$ or non-converge} 
 \State Prediction step : 
    \begin{align*}
        &\pmean_{n+1} = \mean_{n} \qquad \pCov_{n+1} = 2\Cov_n\\
    \end{align*}
\State Generate sigma points :
    \begin{align*}
    &\pp_{n+1}^0 = \pmean_{n+1} \\
    &\pp_{n+1}^j = \pmean_{n+1} + c_j [\sqrt{\pCov_{n+1}}]_j \quad (1\leq j\leq N_\theta)\\ 
    &\pp_{n+1}^{j+N_\theta} = \pmean_{n+1} - c_j [\sqrt{\pCov_{n+1}}]_j\quad (1\leq j\leq N_\theta)
    \end{align*}
\State Analysis step :
   \begin{equation}
   \label{eq:UKI-analysis}
   \begin{split}
        &{\ppy}^j_{n+1} = \G(\bm{\pp}^j_{n+1}) \quad (0\leq j\leq 2N_{\theta})\\ 
        &\widehat{y}_{n+1} = {\ppy}^0_{n+1}\\
         &\pCov^{\theta p}_{n+1} = \sum_{j=1}^{2N_\theta}W_j^{c}
        (\pp^j_{n+1} - \pmean_{n+1} )({\ppy}^j_{n+1} - \widehat{y}_{n+1})^T \\
        &\pCov^{pp}_{n+1} = \sum_{j=1}^{2N_\theta}W_j^{c}
        ({\ppy}^j_{n+1} - \widehat{y}_{n+1} )({\ppy}^j_{n+1} - \widehat{y}_{n+1})^T + 2\Sigma_{\eta}\\
        &\mean_{n+1} = \pmean_{n+1} + \pCov^{\theta p}_{n+1}(\pCov^{pp}_{n+1})^{-1}(y - \widehat{y}_{n+1})\\
        &\Cov_{n+1} = \pCov_{n+1} - \pCov^{\theta p}_{n+1}(\pCov^{pp}_{n+1})^{-1}{\pCov^{\theta p}_{n+1}}{}^{T}\\
    \end{split}
    \end{equation}
    \EndWhile
    \State \Return $\mean_{n+1},\, \Cov_{n+1}$
    \EndFunction
 \end{algorithmic}
\end{algorithm}

\subsection{Theoretical and Numerical Analysis}
\label{ssec:analysis}
In this section, we present theoretical and numerical analysis about the aforementioned calibration method with respect to linear convergence, uncertainty quantification capability, and non-convex optimization. 

\begin{theorem}[Linear Analysis]
\label{TH:LINEAR}
In the context of linear inverse problems, for which $\G(\cdot) = G$, we assume the inverse problem is well-defined, namely $G$ has full column rank, and the initial covariance matrix $\Cov_{0} \succ 0$ is strictly positive definite. 

\begin{itemize}
    \item The iteration for the mean $\mean_n$ and covariance matrix $\Cov_{n}$
converges exponentially fast to limit $\mean_{\infty}, \Cov_{\infty}.$
Furthermore, the limiting mean $\mean_{\infty}$ is a minimizer
of the unregularized least squares
functional $\Phi$~\eqref{eq:LSQ};
the limiting covariance matrix $\Cov_{\infty} = \Big(G^T\Sigma_{\eta}^{-1}G\Big)^{-1}$ is the 
posterior covariance matrix with an uninformative prior.
\item The uncertainty is given as an error bound about the parameter estimation
\begin{equation}
\label{eq:err_bound_linear}
    P\Big( |{\theta_{ref}}_{(i)} - {m_{\infty}}_{(i)}| \leq 3\sqrt{{\Cov_{\infty}}_{(i,i)}} \Big) \geq 99.7\%,
\end{equation}
here the subscript $i$ represents the vector or matrix index, and $\theta_{ref}$ represents the reference parameter, which satisfies $y - G\theta_{ref} \sim \N(0, \Sigma_{\eta})$.
\end{itemize}
\end{theorem}
\begin{proof}
The proof is in \cref{sec:linear-proof}.
\end{proof}

We further study numerically the performance of UKI for nonlinear calibration problems with multiple minimizers. Due to the "discrete" events existing in the EBMs and IBMs, the objective function~(likelihood function) $\Phi(\theta)$ might be lacking of smoothness and might have multiple local optima, which hinders classical optimization approaches. However, we find UKI, as an ensemble-based method, has a better performance in avoiding local optima.

\begin{figure}[ht]
\centering
    \includegraphics[width=0.5\textwidth]{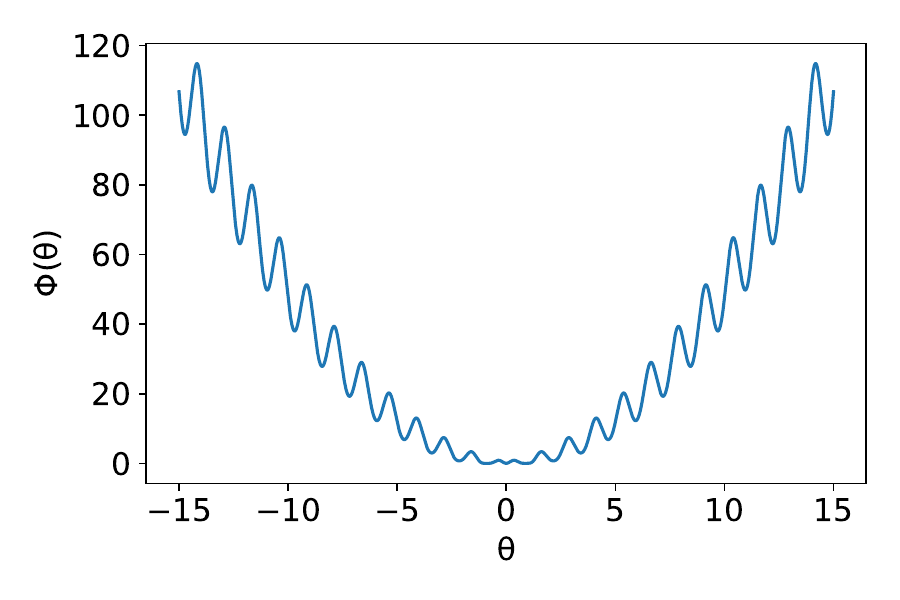}
    \caption{The landscape of the objective function for the illustrative one-dimensional problem.}
    \label{fig:1D-Loss}
\end{figure}

In this comparison, UKI is compared with three other optimization methods, including derivative free ensemble transformed Kalman inversion~(ETKI) described in the Appendix B.2 of the work of Huang \textit{et al.}~\cite{UKI2} and gradient-based Newton and BFGS (Broydon-Fletcher-Goldfarb-Shanno) methods, for a nonconvex one-dimensional nonlinear least square problem~\eqref{eq:LSQ}
with 
$$\G(\theta) = \sin(5\theta) + \theta \qquad y = 0 \qquad \Sigma_{\eta} = \begin{bmatrix}
1.0
\end{bmatrix}.$$
This problem has $5$ global minimizers: $\theta = -0.981,\, -0.821,\, 0,\, 0.821,\, 0.981$, and numerous local minimizers. The landscape of the objective function is depicted in~\cref{fig:1D-Loss}.
Three different initial guesses are chosen:  $\mean_0 = 10$, $\mean_0 = 11$, and $\mean_0 = 12$, and for UKI and ETKI the initial covariance is $\Cov_0 = \I$. The ensemble size of the ETKI is $J=10$, which is larger than the number of $\sigma$-points used in UKI. For Newton and BFGS methods, we use library Optim.jl~\cite{mogensen2018optim} with the line search method proposed in the work of Hager $\&$ Zhang~\cite{hager2006algorithm}. 

\begin{figure}[ht]
\centering
    \includegraphics[width=0.75\textwidth]{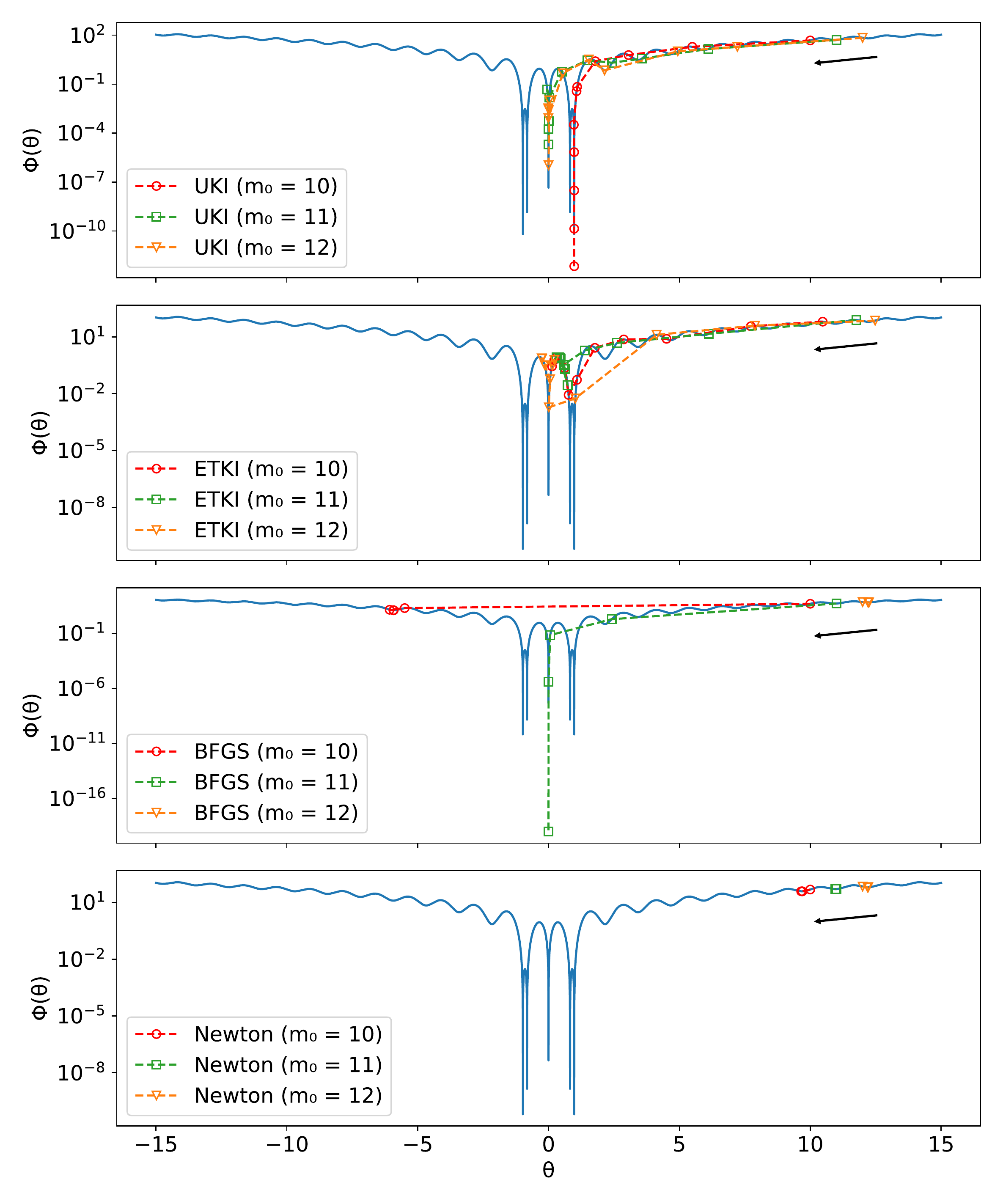}
    \caption{Convergence histories of different optimization approaches initialized at $m_0 = 10$ (red circle dashed lines), $m_0 = 11$ (green square dashed lines), and $m_0 = 12$ (orange triangle dashed lines) for the illustrative one-dimensional problem.}
    \label{fig:1D-Opt}
\end{figure}

The convergence histories obtained with these approaches are depicted in~\cref{fig:1D-Opt}. Both Newton and BFGS methods are prone to be trapped in local optima, primarily because only local gradient and Hessian information are used in the solution process. 
However, neither ETKI nor UKI are trapped in local optima. For UKI, although there is no guarantee that the global minimizer can be found, the ensemble $\{\theta_{n+1}^{j}\}$ around the mean value $\pmean_{n+1}$ brings non-local information and has averaging effect~\cite{UKI1} on the landscape. And therefore, it is able to avoid local optima. Compared to UKI, ETKI suffers from random noise, and the estimated means oscillate around the optima.

\section{Fluid Structure Interaction under Embedded Boundary Framework}
\label{sec:FSI}
\subsection{Governing Equations}
We consider solving fluid-structure interaction problems involving large structural deformations and even structural damage. Let $\Omega_{F}$ and $\Omega_{S}$ denote the fluid and the structural subdomains. The fluid is assumed to be compressible and inviscid, governed by the following Euler equations in the Eulerian frame
\begin{equation}
    \frac{\partial W(x,t)}{\partial t} + \nabla \cdot \mathcal{F}(W) = 0, \qquad \forall x \in \Omega_{F}(t), \, t>0,
\label{eq:EulerEquation}
\end{equation}
\begin{equation*}
W=\begin{bmatrix} \rho \\ \rho v \\ \rho e_{t} \end{bmatrix}, \qquad \mathcal{F} = \begin{bmatrix} \rho v^T \\  \rho v \otimes v + p \I \\ (\rho e_{t} + p)v^T \end{bmatrix}, 
\end{equation*}
where $\rho$, $v$, $p$ and $e_{t}$ denote the fluid density, velocity, pressure, and the total energy per unit mass, respectively. $\I$ denotes an identity matrix. The fluid is assumed to be a perfect gas, and the following equation of state is applied to close~\cref{eq:EulerEquation}
\begin{equation}
p=\left( \gamma_{G} -1 \right) \rho e ,
\label{eq:air_EOS}
\end{equation}
where the heat capacity ratio $\gamma_{G}$ is set to be $1.4$ in this work, and $e$ denotes the internal energy per unit mass. 

The dynamics of structure is governed by the following equation of motion in the Lagrangian frame
\begin{equation}
\rho_{s} \ddot{u}(X,t) - \nabla \cdot \sigma = b, \qquad \forall X \in \Omega_{S}(0), \, t > 0, \label{eq:SolidEquation}
\end{equation}
where $u$ denotes displacement, $\rho_{s}$ the material's mass density and  $\sigma$ the Cauchy or Piola-Kirchhoff stress tensor. $b$ denotes the body force acting in $\Omega_{S}$, which is assumed to be zero in this study. The dot above a variable represents partial derivative with respect to time. Given a structural material of interest, the closure of~\eqref{eq:SolidEquation} is performed by specifying
a constitutive law that relates the stress tensor $\sigma$ to the
strain tensor~$\epsilon$.

% We consider that the structure is linear elastic and isotropic, and the constitutive equation is given by 
% \begin{equation}
% \sigma = \frac{E\nu}{(1+\nu)(1-2\nu)} tr(\mathbf{\bm{\epsilon}})\mathbf{I} + \frac{E}{1+\nu} \bm{\epsilon},
% \end{equation}
% where $E$ and $\nu$ are the Young's modulus and Poisson's ratio of the material, respectively. $tr(\cdot)$ denotes the trace operator. 

% The structural damage is modeled using a continuum damage mechanics model. Specifically, a scalar damage state variable, $\omega(\mathbf{X},t)$, is introduced to represent the damage. The damage effect is imposed by modifying the elastic moduli ($E^{d}$ and $G^{d}$) as follows:
% \begin{eqnarray}
%     E^{d}(\mathbf{X},t) &=& E_{0}\left( 1-\omega(\mathbf{X},t) \right),\\
%     G^{d}(\mathbf{X},t) &=& G_{0}\left( 1-\omega(\mathbf{X},t) \right),
% \end{eqnarray}
% where $E_{0}$ and $G_{0}$ are the Young's modulus and shear modulus of the material without damage. In this work, we assume $\omega(\mathbf{X},t)$ varies from $\omega_{min} = -0.1$ (hardening) to $\omega_{max} = 0.9$ (softening).

The fluid-structure interface, $\Gamma_{FS} = \partial\Omega_{F}(t) \cap \partial \Omega_{S}(t)$, is assumed to be impermeable, at which we enforce the continuity of normal velocity and surface traction, i.e.,
\begin{equation}
\begin{split}
\label{eq:IC2}
(v - \dot{u})\cdot n = 0, \qquad &\text{on} \quad \Gamma_{FS}, \\ 
-pn = \sigma \cdot n, 
\qquad &\text{on} \quad \Gamma_{FS},
\end{split}
\end{equation}
where $n$ denotes the unit normal to $\Gamma_{FS}$. 

\subsection{An Embedded Boundary Computational Framework}
\label{subsec:computational_framework}
The above coupled problem is solved by a recently developed fluid-structure coupled computational framework\cite{farhat2010robust,wang2011algorithms,wang2012computational,main2017enhanced,huang2018family}---AERO-Suite~\footnote{\url{https://bitbucket.org/frg/}}. The framework couples a finite volume CFD solver with a finite element computational structural dynamics solver using an embedded boundary method and a partitioned coupling procedure. The framework has been verified and validated on various applications, including dynamic implosion of underwater cylindrical shells~\cite{wang2011algorithms,wang2012computational,farhat2013dynamic}, F/A-18 vertical tail buffeting~\cite{lakshminarayan2014embedded}, and Mars landing parachute inflation dynamics~\cite{huang2020embedded,huang2020modeling,huang2021homogenized}. 

\subsubsection{FIVER: A Finite Volume Method Based on Exact Riemann Solver}
\begin{figure}[ht]
\centering
    \includegraphics[width=0.8\textwidth]{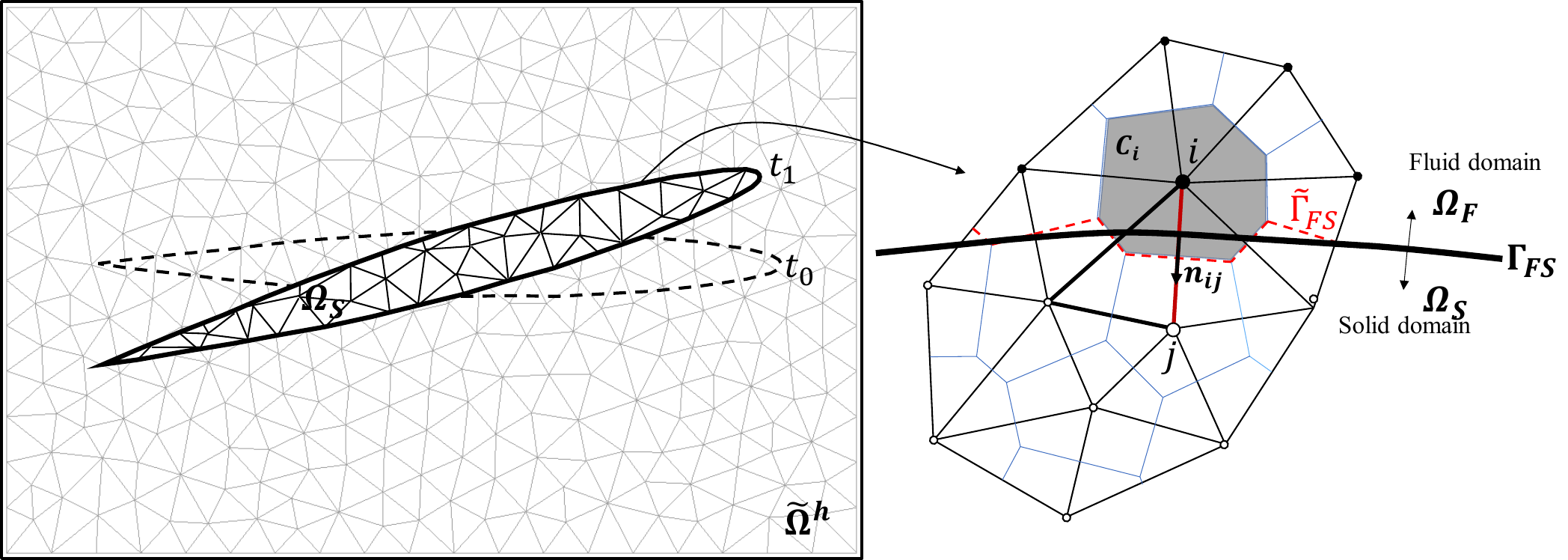}
    \caption{Illustration of the non-interface-conforming mesh under embedded boundary framework.}
    \label{fig:EBM_mesh}
\end{figure}

The fluid governing equation is semi-discretized in an augmented fluid domain $\widetilde{\Omega} = \Omega_F \cup \Omega_S$, using an unstructured, node centered, non-interface-conforming finite volume mesh, denoted by $\widetilde{\Omega}^{h}$. Figure~\ref{fig:EBM_mesh} presents an example to illustrate the non-interface-conforming mesh. Within each control volume ($\Cv_i$) in $\widetilde{\Omega}^{h}$, \cref{eq:EulerEquation} is integrated as
\begin{equation}
\frac{\partial W_i}{\partial t} + \frac{1}{\lVert \Cv_i \rVert} \sum\limits_{j\in Nei(i)} \int_{\partial \Cv_{ij}} \mathcal{F}(W)\cdot n_{ij} dS = 0,
\label{eq:discrete_Euler}
\end{equation}
where $W_i$ denotes the average of $W$ in $\Cv_i$, $\lVert \Cv_i \rVert$ denotes the volume of $\Cv_i$,  $Nei(i)$ denotes the set of nodes connected to node $i$ by an edge, $\partial \Cv_{ij} = \partial \Cv_i \cap \partial \Cv_j$, and $n_{ij}$ is the unit normal to $\partial \Cv_{ij}$. Away from the fluid-structure interface, the flux through each segment $\partial \mathcal{\Cv}_{ij}$ in~\cref{eq:discrete_Euler} is approximated using Roe’s (or any
other similar) approximate Riemann solver~\cite{roe1981approximate} equipped with a MUSCL technique~\cite{van1979towards}
and a slope limiter.

At the embedded fluid-structure interface, the flux in~\cref{eq:discrete_Euler} is approximated with the reconstructed fluid state vector, which satisfies the transmission conditions~\eqref{eq:IC2}. The reconstruction is through solving a one-dimensional fluid-solid Riemann problem.
Specifically, as shown in the right figure in Figure~\ref{fig:EBM_mesh}, suppose node $i$ belongs to the fluid subdomain, node $j$ belongs to the solid subdomain, and edge $i$-$j$ intersects the embedded interface $\Gamma_{FS}$. 
We construct the following one-dimensional fluid-solid Riemann problem along $n_S$, the unit normal of $\Gamma_{FS}$ (toward the structure) at its intersection with edge $i$-$j$, at each time-instance $t_n$ 

\begin{eqnarray}
\label{eq:FSRiemann}
\frac{\partial {w}}{\partial \tau} + \frac{\partial {\mathcal{F}}({w})}{\partial \xi} &=& 0, \quad\quad 0\leq\tau\leq \Delta t^n,~\xi<(\dot{u}_0\cdot n_S)\tau, \\
{w}(\xi,0) &=& w_i, \quad\quad \mbox{$\xi<0$,}\\
v((\dot{u}_0\cdot n_{S})\tau,\tau) &=& \dot{u}_0 \cdot  n_S, \quad\quad 0\leq\tau\leq \Delta t_n, 
\end{eqnarray}
$\xi$ is the spatial coordinate along the one-dimensional axis aligned with $ n_{S}$ and centered at the intersection point. The initial state $w_i$ is the projection of $W_i$ on $ n_S$, i.e.
\begin{equation}
w_i = 
\begin{bmatrix}
\rho_i \\
\rho_i (v_i \cdot  n_{S}) \\
\rho_i \big(e_i + \frac{1}{2} (v_i \cdot  n_{S})^2\big)\\
\end{bmatrix}.
\end{equation}
$\dot{u}_0$ denotes the velocity of the structure at $t_n$. The exact solution of this Riemann problem can be derived analytically. And the state variables and the tangential fluid velocity are used to reconstruct the fluid state vector at the fluid-structure interface. 
The resulting semi-discretization of \cref{eq:EulerEquation} can be written in a compact form as
\begin{equation}
\frac{d \mathbf{W}^h}{dt} + \mathbf{V}^{-1}\mathbf{F}(\mathbf{W}^h) = 0,
\label{eq:discrete_Euler_final}
\end{equation}
where $\mathbf{W}^h$, $\mathbf{V}$, and $\mathbf{F}(\mathbf{W}^h)$ denote the vector of semidiscrete fluid state variable, the diagonal matrix storing the volume of control volumes, and the vector of numerical flux, respectively.

\subsubsection{A Finite Element Structural Solver}
The structural governing equation is solved using a standard Galerkin finite element method. The semi-discretized equation is written as
\begin{equation}
\mathbf{M}\frac{\partial^2 \mathbf{u}^h}{\partial t^2} + \mathbf{f}^{int}\Big(\mathbf{u}^h, \frac{\partial \mathbf{u}^h}{\partial t}\Big) = \mathbf{f}^{ext},
\label{eq:discrete_Solid}
\end{equation}
where $\mathbf{M}$ denotes the mass matrix, $\mathbf{u}^h$ denotes the discrete displacement vector. $\mathbf{f}^{int}$ and $\mathbf{f}^{ext}$denote the discrete internal force and external force vector, respectively. Based on the dynamics interface condition (\cref{eq:IC2}), the fluid-induced external forces are computed by integrating the fluid pressure over the embedded interface. Additional details can be found in Section 3.8.3 in the work of Main~\cite{main2014implicit}.

\subsubsection{Partitioned Coupling Procedure}
We employ a partitioned procedure to advance the fluid subsystem~\eqref{eq:discrete_Euler_final} and structure subsystem~\eqref{eq:discrete_Solid} in time, following the work of Farhat~\textit{et al.}~\cite{farhat2010robust}. 
Specifically, the fluid and structural discretized equations are solved independently with different time integrators.
The two solvers exchange information at the fluid-structure interface once per time step in a staggered manner, i.e., the fluid and solid time steps are offset by half a step. This is a designed feature to achieve second-order accuracy in time while maintaining optimal numerical stability.

\section{Numerical Results}
\label{sec:app}
In this section, we apply the unscented Kalman inversion and the embedded boundary framework presented in Sections~\ref{sec:Bayesian} and \ref{sec:FSI} to calibrate parameters, as well as quantify their uncertainty, for two FSI problems with noisy observation data. The first case is a one-dimensional piston model problem, where the damping coefficient, the spring stiffness, and the initial fluid pressure are calibrated from the piston displacement. The second case is a challenging three-dimensional aeroelasticity problem of a damaged aircraft wing, where the damage field, specifically the coefficients of first 5 modes, is inferred from the displacement during its transonic buffet. For both cases, we prescribe a reference value for the parameters of interest and generate the time-series observation data in which additional Gaussian random noise is added. Then, using the noisy data, we solve the 
Bayesian calibration problem to retrieve the predefined parameters and associated uncertainties.
The code associated with these two test cases is accessible online:
\url{https://github.com/Zhengyu-Huang/InverseProblems.jl}.

\subsection{Piston Problem}
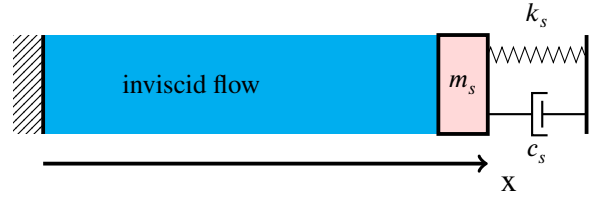
\begin{figure}[ht]
\centering
    \begin{tikzpicture}[scale=1.3]
    \tikzstyle{damper}=[thick,decoration={markings,  
   mark connection node=dmp,
   mark=at position 0.5 with 
   {
     \node (dmp) [thick,inner sep=0pt,transform shape,rotate=-90,minimum
 width=15pt,minimum height=3pt,draw=none] {};
     \draw [thick] ($(dmp.north east)+(2pt,0)$) -- (dmp.south east) -- (dmp.south
 west) -- ($(dmp.north west)+(2pt,0)$);
     \draw [thick] ($(dmp.north)+(0,-5pt)$) -- ($(dmp.north)+(0,5pt)$);
   }
 }, decorate]
      \fill [pattern = north east lines] (-0.3,0) rectangle (0.3,1.0);
      \fill[cyan] (0,0) rectangle (4,1.0);
      \draw [fill=red!15!white,line width=0.5mm] (4,0) rectangle (4.5,1);
      \draw node at (1.5,0.5) {inviscid flow};
      \draw[decoration={aspect=0.3, segment length=1.5mm, amplitude=1mm,zigzag},decorate] (4.5,0.8) -- (5.5,0.8);
      \draw[damper] (4.5,0.2) -- (5.5,0.2);
    %   \fill [pattern = north east lines] (5.5,0) rectangle (6.0,1.0);
      \draw[line width=0.5mm](0, 0)-- (0, 1);
      \draw[line width=0.5mm](5.5, 0)-- (5.5, 1);
      \draw[thick,->,line width=0.5mm,](0,-0.3) -- (4.5,-0.3) node[anchor=north west][scale=1.2] {x};
      \draw node at (4.25,0.5) {$m_s$};
      \draw node at (5.0,1.2) {$k_s$};
      \draw node at (5.0,-0.2) {$c_s$};
    \end{tikzpicture}
    \caption{One-dimensional piston system}
    \label{fig:piston-system}
\end{figure}

We first consider a canonical FSI model problem: a one-dimensional piston problem depicted in \cref{fig:piston-system}. 
The inviscid fluid is governed by the one-dimensional Euler equations defined within domain $x \in [0, 1 - u]$, where $u$ is the displacement of the piston.
The fluid is initially at rest, with initial density, velocity, and pressure given by:
\begin{equation*}
    \rho_0 = 1.225, \qquad v_0 = 0, \qquad p_0 = 2.0.
\end{equation*}
No-penetration wall boundary condition is imposed on the left end of the fluid domain.
On the right, the fluid interacts with a piston. The piston is connected with a spring and a damper that are attached to a wall on the right at $x=1$. The dynamics of the piston can be described by the second order ordinary differential equation:
\begin{equation*}
    m_s \ddot{u} + c_s \dot{u} + k_s {u} = f^{ext},
\end{equation*}
where the mass coefficient is $m_s=1$, the damping coefficient is $c_s = 0.5$, and the spring stiffness is $k_s = 2$. The external force $f_{ext}$ is the fluid pressure load on the piston:
$$ f^{ext}(t) = p(x_p(t)),$$
where $x_p(t) = 1 - u(t)$.
The piston is initialized at $x=1$, with $u(0) = 0$ and $\dot{u}(0) = 0$.

For the forward problem, the computational domain is semi-discretized by a uniform grid with $\Delta x = 5\times10^{-3}$. The explicit 2nd-order Runge-Kutta time integrator~\cite{} is used for the fluid and implicit 2nd-order mid-point rule time integrator is applied for the structure with a constant time step $\Delta t = 10^{-3}$. The coupling between fluid and structure is based on a partitioned procedure and an embedded boundary method which are applied in FIVER (Section~\ref{subsec:computational_framework}). The coupled system is integrated till the final time $T = 1$.

\begin{figure}[ht]
\centering
    \includegraphics[width=0.45\textwidth]{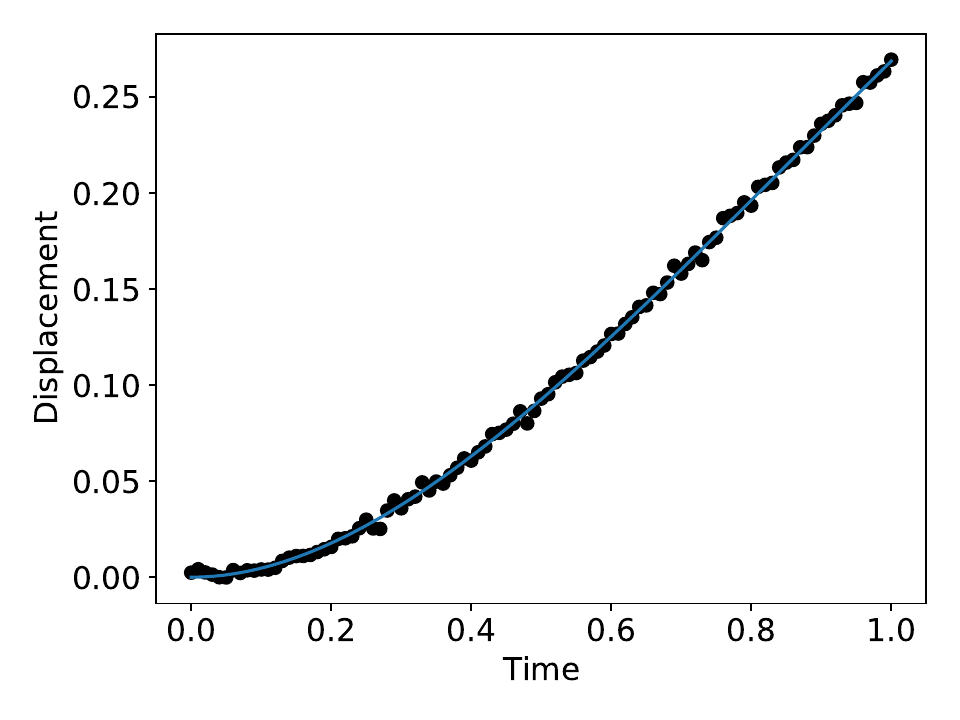}
    \includegraphics[width=0.45\textwidth]{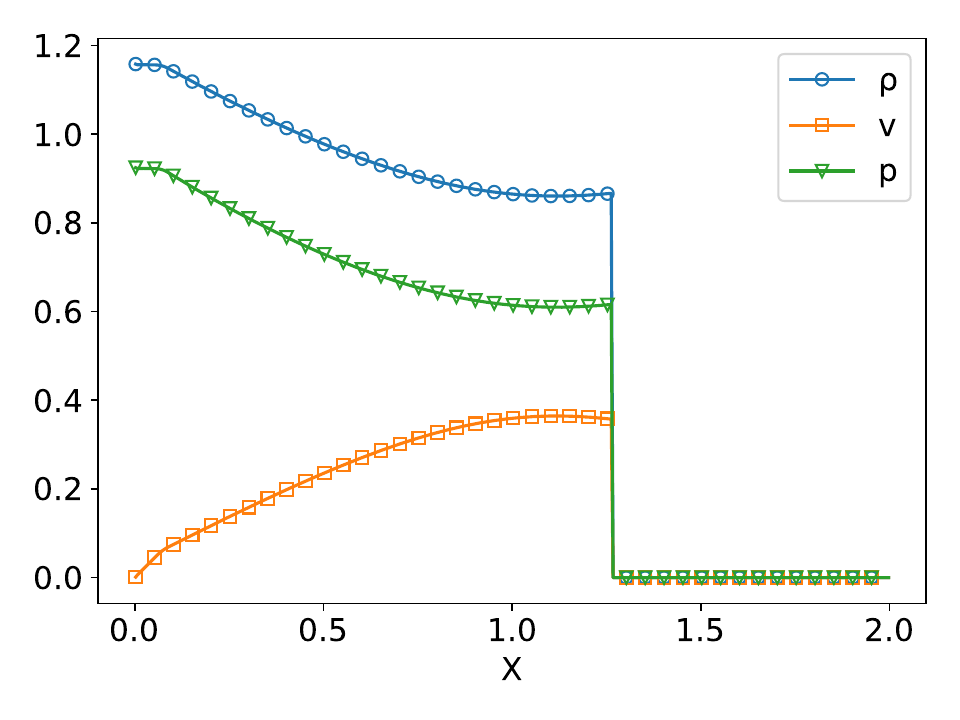}
    \caption{Left: the displacement of the piston with $100$ noisy measurements~(black dots); Right: the fluid states at the end time.}
    \label{fig:piston-forward}
\end{figure}

For the inverse problem, we consider two scenarios:
\begin{enumerate}
    \item calibrate only the structure parameters $\theta = [c_s;\,k_s]$ with $N_{\theta} = 2$;
    \item calibrate both structure parameters and fluid initial pressure $\theta = [c_s;\,k_s;\,p_0]$ with $N_{\theta} = 3$.
\end{enumerate}
The observation data consist of $N_y=100$ piston displacement measurements collected at every $10^{-2}$ till $T = 1$. These observation data are generated with reference parameters and corrupted with zero-mean random Gaussian noise with standard deviation $\sigma_{\eta} = 2\times10^{-3}$. The observation data are depicted in~\cref{fig:piston-forward}-left; and the corresponding fluid states at the end time is depicted in~\cref{fig:piston-forward}-right , the rarefaction wave is generated in the flow due to the receding motion of the piston.

\begin{figure}[ht]
\centering
    \includegraphics[width=0.45\textwidth]{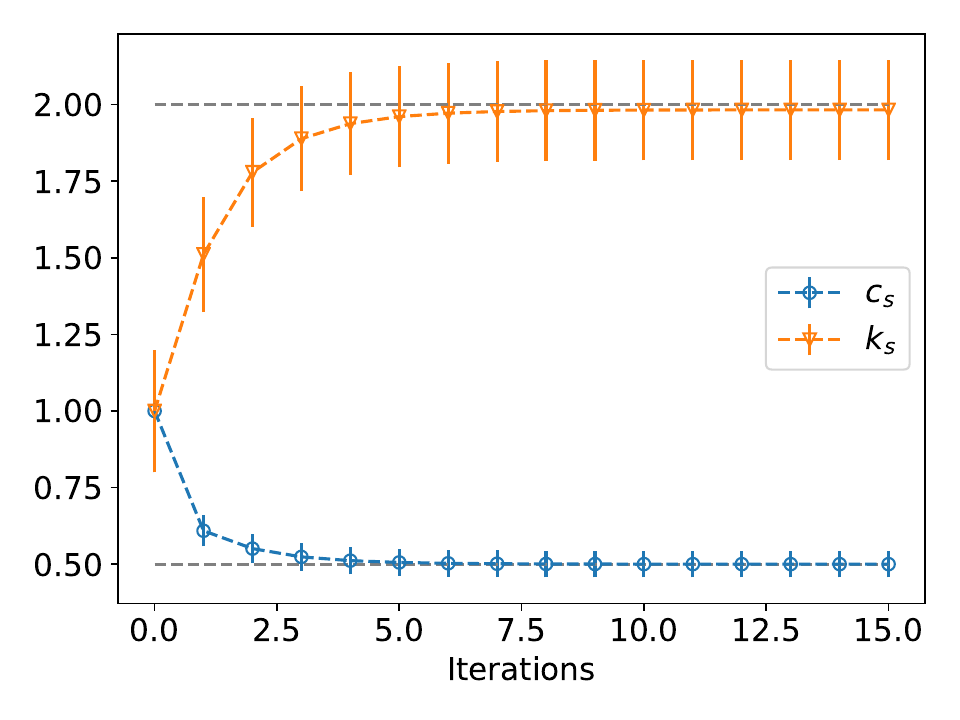}
    \includegraphics[width=0.45\textwidth]{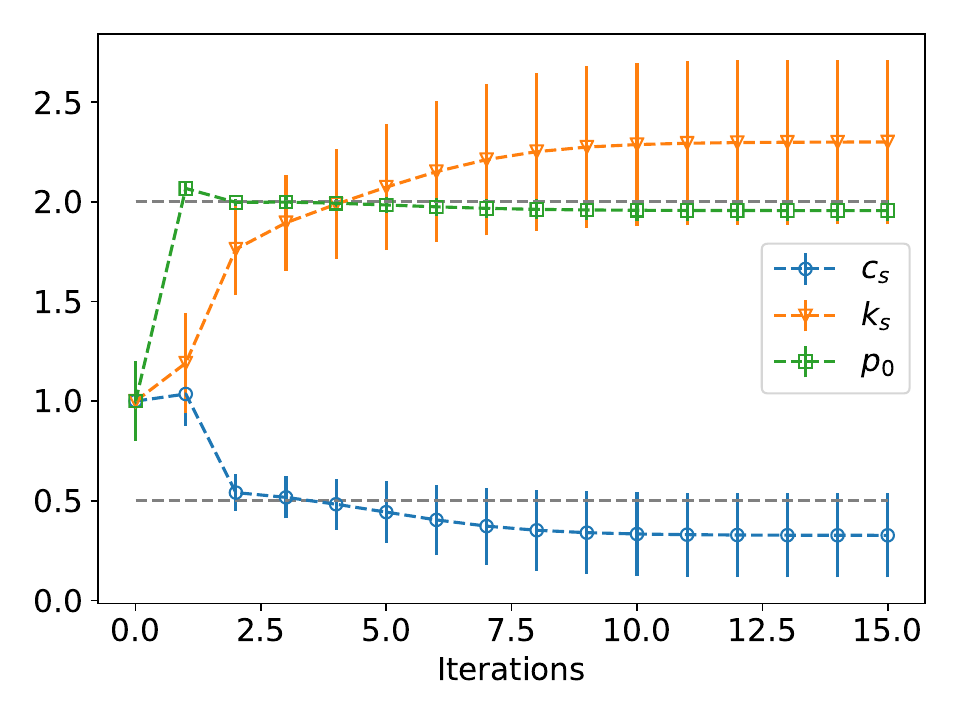}
    \caption{Convergence of the piston model problem with UKI for both 2-parameter~(left) and 3-parameter~(left) scenarios; the true parameter values are represented by dashed grey lines. The error bar associated with each data point denotes the $2-\sigma$ confidence intervals.}
    \label{fig:piston-uki}
\end{figure}

\begin{figure}[ht]
\centering
    \includegraphics[width=0.75\textwidth]{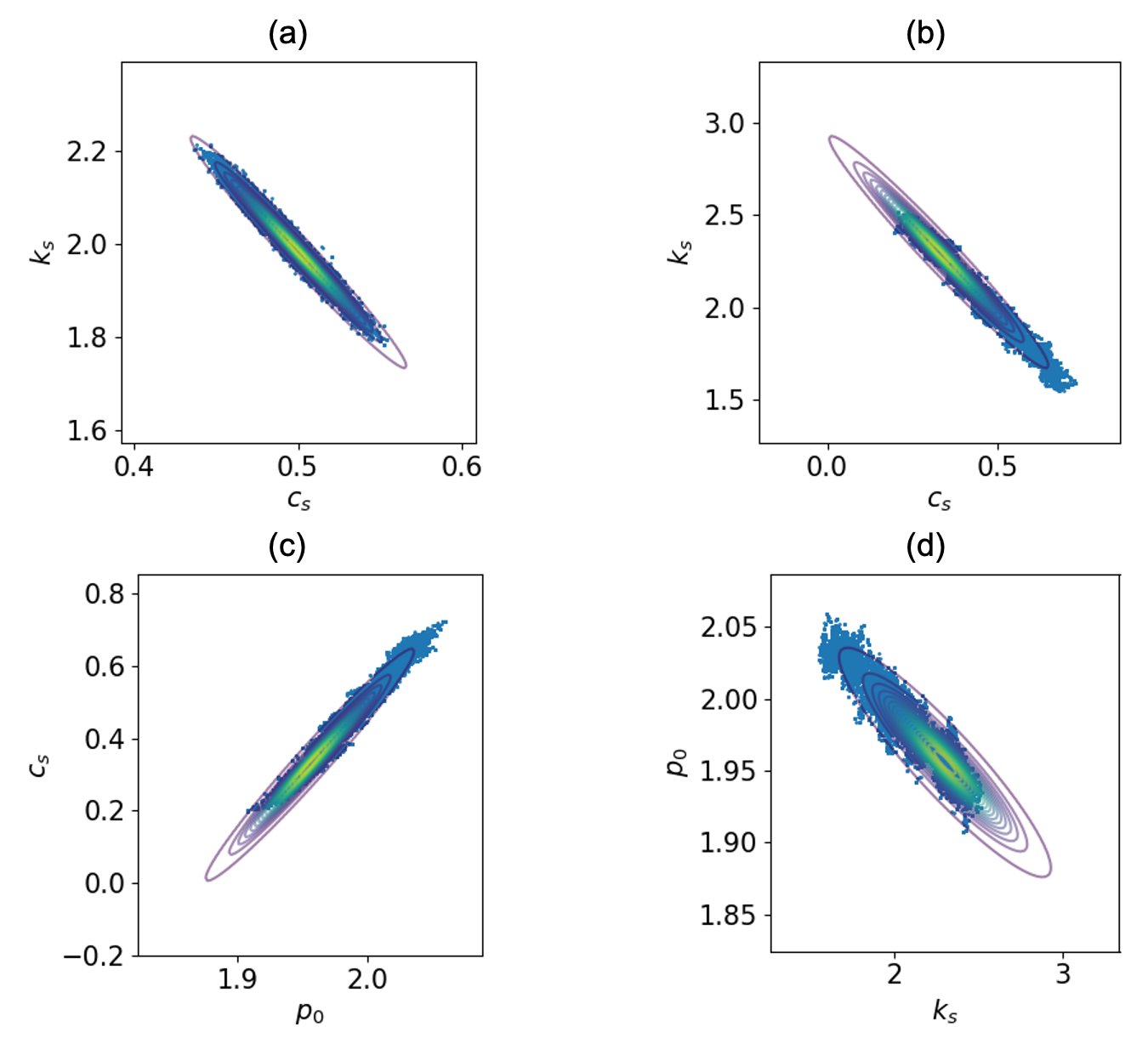}
    \caption{Pairwise posterior distribution of the piston model problem for both 2-parameter~(a) and 3-parameter~(b-d) scenarios. Contour plot: the posterior distribution obtained by the UKI at the 15th iteration; blue dots: reference posterior distribution obtained by MCMC.}
    \label{fig:piston-uki-mcmc}
\end{figure}

For both scenarios, the UKI is initialized with $\theta_0 \sim \N(\mathds{1}, 0.1^2\I)$. 
The estimated parameters and associated 2-$\sigma$ confidence intervals for each component at each iteration are depicted in \cref{fig:piston-uki}. 
The UKI converges efficiently.  Reference values falls within the confidence interval with high probability, although the uncertainty is larger for the 3-parameter scenario.
% The estimation of the parameters at the 10th iteration is 
% \begin{equation*}
%     \begin{bmatrix}
%     c_s\\
%     k_s
%     \end{bmatrix}\sim \N
%     \Big(
% \begin{bmatrix}
% 0.499\\ 
% 1.97
% \end{bmatrix},
% \begin{bmatrix}
% 0.00156 &-0.00593\\  
% -0.00593 & 0.0232
% \end{bmatrix}
% \Big).
% \end{equation*}
%
The reference posterior distribution is approximated by the random walk MCMC method with a step size $1.0\times10^{-2}$ and  $5\times10^4$ samples (with a $10^4$ sample burn-in period). 
The posterior distributions obtained by the UKI at the 15th iteration are depicted in~\cref{fig:piston-uki-mcmc}. 
The UKI delivers a very similar posterior distribution, but at a significantly cheaper computational cost.

% \begin{equation*}
%     \begin{bmatrix}
%     c_s\\
%     k_s\\
%     p_0
%     \end{bmatrix}\sim \N
%     \Big(
% \begin{bmatrix}
% 0.32666738891228536\\
% 2.30011362238824\\
% 1.9558503327585004\\ 
% \end{bmatrix},
% \begin{bmatrix}
% 0.011164118490875523 &-0.02137571156507171& 0.0027199429843845413 \\ -0.021375711565071717 &0.04268101835345774& -0.004999749253725635\\ 0.0027199429843845417 &-0.004999749253725636& 0.0006905475898639387
% \end{bmatrix}
% \Big).
% \end{equation*}

\subsection{Wing Damage Detection Problem}
\label{ssec:damage}
\begin{figure}[ht]
\centering
    \includegraphics[width=0.8\textwidth]{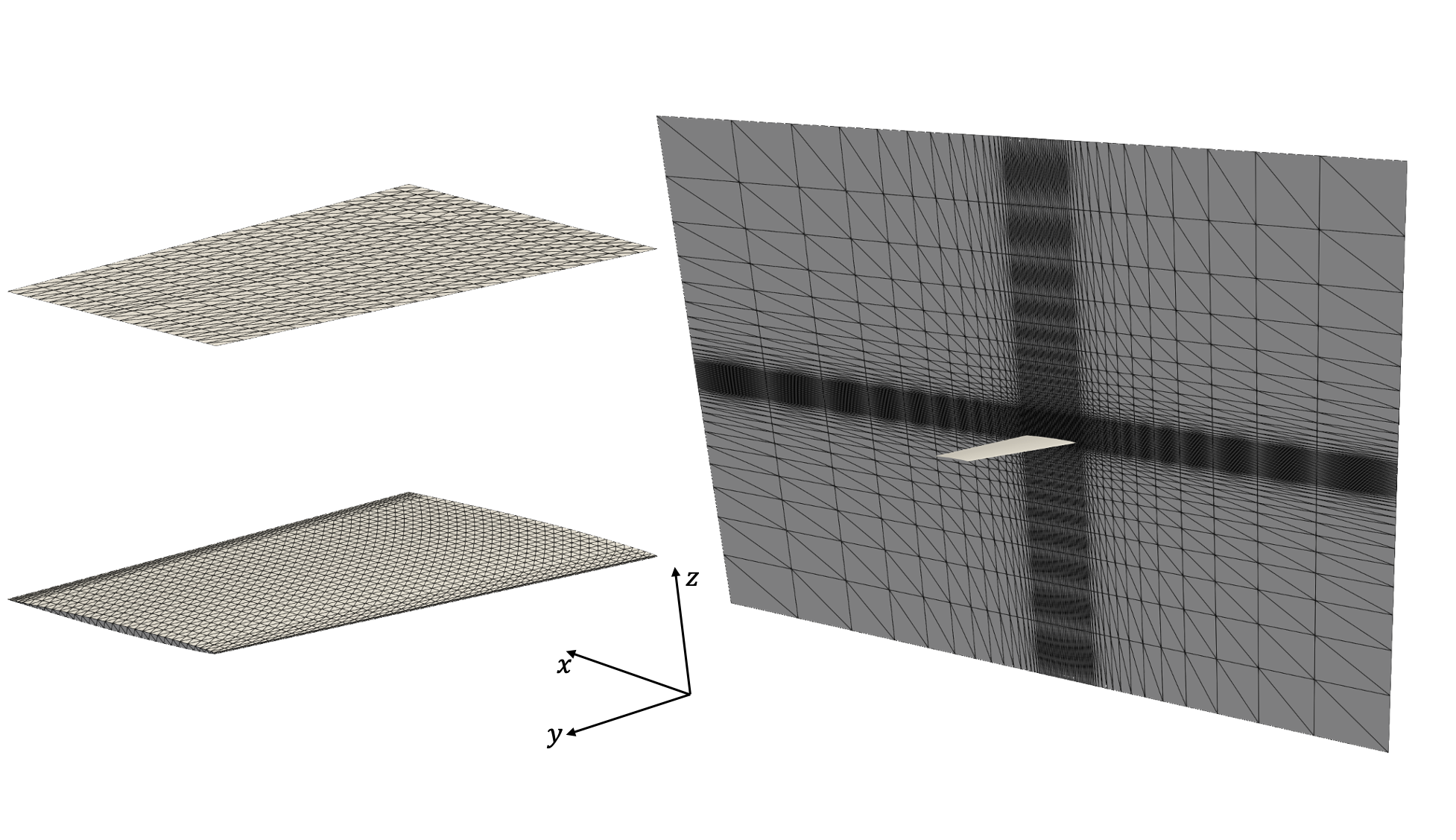}
    \caption{Computational models for the AGARD wing: finite element structural model~(left-top), the wing skin surface of this model (left-bottom), and CFD mesh (right).}
    \label{fig:Wing-mesh}
\end{figure}

The second test is a challenging real-world FSI problem associated with a damaged AGARD wing~\cite{yates1987agard} undergoing transonic buffet. We consider a cruise condition, where the atmospheric density is 
$\rho_{\infty} = 0.61\times10^{-7}~\textrm{slug/in}^3$ 
and the atmospheric pressure is 
$p_{\infty} = 6~\textrm{slugs}/(\textrm{in} \cdot \textrm{s}^2)$; the free-stream Mach number is $M_{\infty} = 0.97$; and the angle of attack is $\alpha = 7.5^\circ$.
The AGARD wing structure is modeled as a nonlinear elastic composite shell~(See~\cref{fig:Wing-mesh}-left). The orthotropic properties of this material are density $\rho_s = 4.6\times10^{-4}~\textrm{slug/in}^3$, parallel Young's modulus $E_1 = 5.874\times10^6~\textrm{slug}/(\textrm{in}\cdot \textrm{s}^2)$, orthogonal Young's modulus 
$E_2 = 7.20\times10^6~\textrm{slug}/(\textrm{in}\cdot \textrm{s}^2)$,
rigidity modulus $G = 7.1688\times10^5  \textrm{slug}/(\textrm{in}\cdot \textrm{s}^2)$, 
and Poisson's ratio $\nu = 0.31$.
The physical dimensions are reported in \cref{ta:specAGARD}. 
The thickness distribution is governed by the airfoil shape. 
\begin{table}[!ht]
        \caption{Geometrical properties of the AGARD wing.}
        \centering
        \begin{tabular}{l c}
         \hline\Xhline{2\arrayrulewidth}
          Parameter                &  Type/Value      \\ [0.5ex]
          %heading
          \hline\Xhline{2\arrayrulewidth}
          \hspace{5pt}Wingspan       &   30 inches      \\ [1ex]
          \hspace{5pt}Root chord    &   21.96  inches  \\ [1ex]
          \hspace{5pt}Tip chord     &   14.5  inches   \\ [1ex]
          \hspace{5pt}Sweep          &   45$^{\circ}$   \\ [1ex]
          \Xhline{2\arrayrulewidth}
        \end{tabular}
      \label{ta:specAGARD}
\end{table}

The computational fluid domain $[-100'',\,150'']\times[0'',\,150'']\times[-100'',\,100'']$ is discretized by
the three-dimensional tetrahedral Eulerian mesh. The mesh contains $218,880$ vertices and $1,248,912$ tetrahedral elements, with a $0.75$ inch resolution near the AGARD wing~(See~\cref{fig:Wing-mesh}-right). 
The AGARD wing finite element model~\cite{lesoinne2001linearized} consists of $800$ triangular composite shell elements and $2,646$ degrees of freedom~(See~\cref{fig:Wing-mesh}-left-top). The embedded surface, representing the wing skin surface, consists of 2788 triangles~(See~\cref{fig:Wing-mesh}-left-bottom).

\begin{figure}[ht]
\centering
    \includegraphics[width=0.6\textwidth]{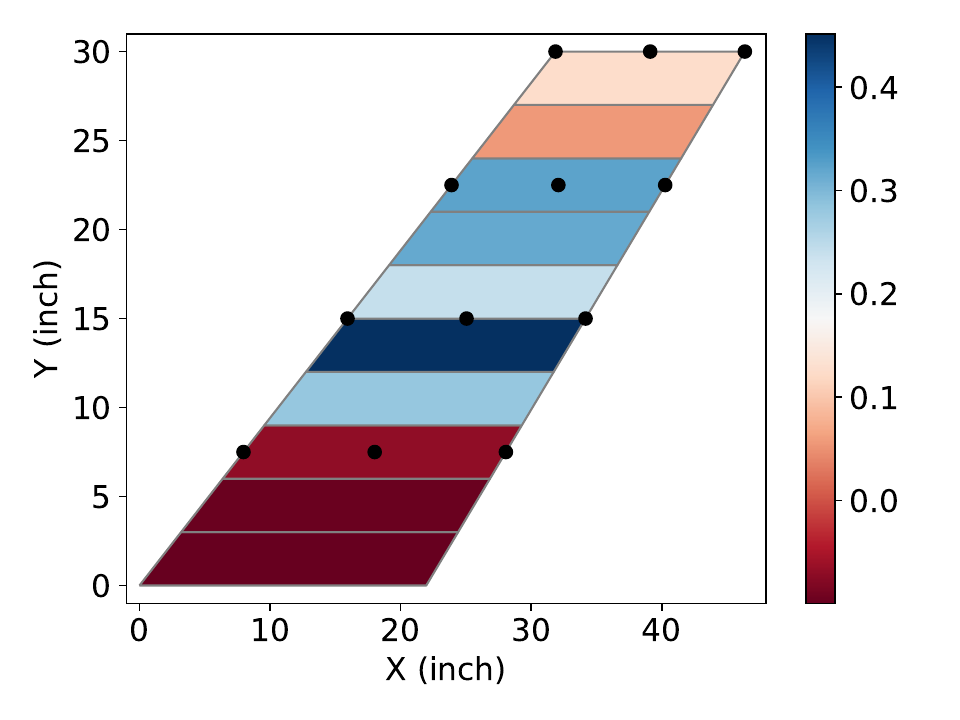}
    \caption{Reference damage field $\omega_{ref}(y)$ of the AGARD wing and the 12 pointwise measurement locations (black dots).}
    \label{fig:Wing-Damage-Obs}
\end{figure}

We assume the damage in aircraft wing to be isotropic elasticity-based damage and consider only the damage in the spar direction~($y$-direction) for simplicity. Specifically, a simple continuum damage model is used in which a scalar damage variable $\omega$ is used to measure the average damage effects on the material's mechanical response. The damage effects are imposed by modifying the elastic moduli of the material:
\begin{equation*}
    E_1^d(y, \theta) = \left(1 - \omega (y, \theta)\right) E_{1}\qquad 
    E_2^d(y, \theta) = \left(1 - \omega (y, \theta)\right) E_{2}\qquad \textrm{and} \qquad
    G^d(y, \theta) = \left(1 - \omega (y, \theta)\right) G.
\end{equation*}
Here, the damage $\omega (y, \theta)$ is designed to vary between $\omega_{min} = -0.1$ (hardening) to $\omega_{max} = 0.9$ (softening):
\begin{equation}
\label{eq:omega-damage}
    \omega(y, \theta) = \frac{\omega_{max} - \omega_{min}}{1 - \frac{\omega_{max}}{\omega_{min}} \exp\left( \log a(y, \theta)\right)} + \omega_{min},
\end{equation}
where $\log a(y, \theta)$ is a log-Gaussian random field that depends on parameters $\theta\in\R^{N_\theta}$. Specifically, the log-Gaussian random field is approximated by the following Karhunen-Lo\`{e}ve~(KL) expansion
\begin{equation}
\label{eq:KL-1d}
    \log a(y,\theta) = \sum_{l=1}^{+\infty} \theta_{(l)}\sqrt{\lambda_l} \psi_l(y),
\end{equation}
where the random variables $\theta_{(l)} \sim \N(0,1)$ are independent and identically distributed, and the eigenpairs are of the form
\begin{equation*}
    \psi_l(y) = \sqrt{2}\cos(\pi l y) \quad \textrm{ and }
                 \quad \lambda_l = (\pi^2 l^2 + \tau^2)^{-d} \quad \text{with} \quad \tau=2, \quad d=1.
\end{equation*}
%Here, we consider the inverse length scale and the regularity of the random field to be $\tau=2$ and $d=1$. 
This setup gives rise to a zero mean state of $\log a(y, \theta)$ which corresponds to the undamaged state $\omega(y, \theta) = 0$. And the damage $\omega$ increases monotonically with $\log a$. In this test, we consider the first 10 KL modes (i.e., with 10 sampled $\theta_{(l)}$) and generate the reference damage field $\omega_{ref}(y)$ from \cref{eq:omega-damage,eq:KL-1d}. The resulting damage field $\omega_{ref}(y)$ is depicted in \cref{fig:Wing-Damage-Obs}. 

\begin{figure}[ht]
\centering
    \includegraphics[width=0.45\textwidth]{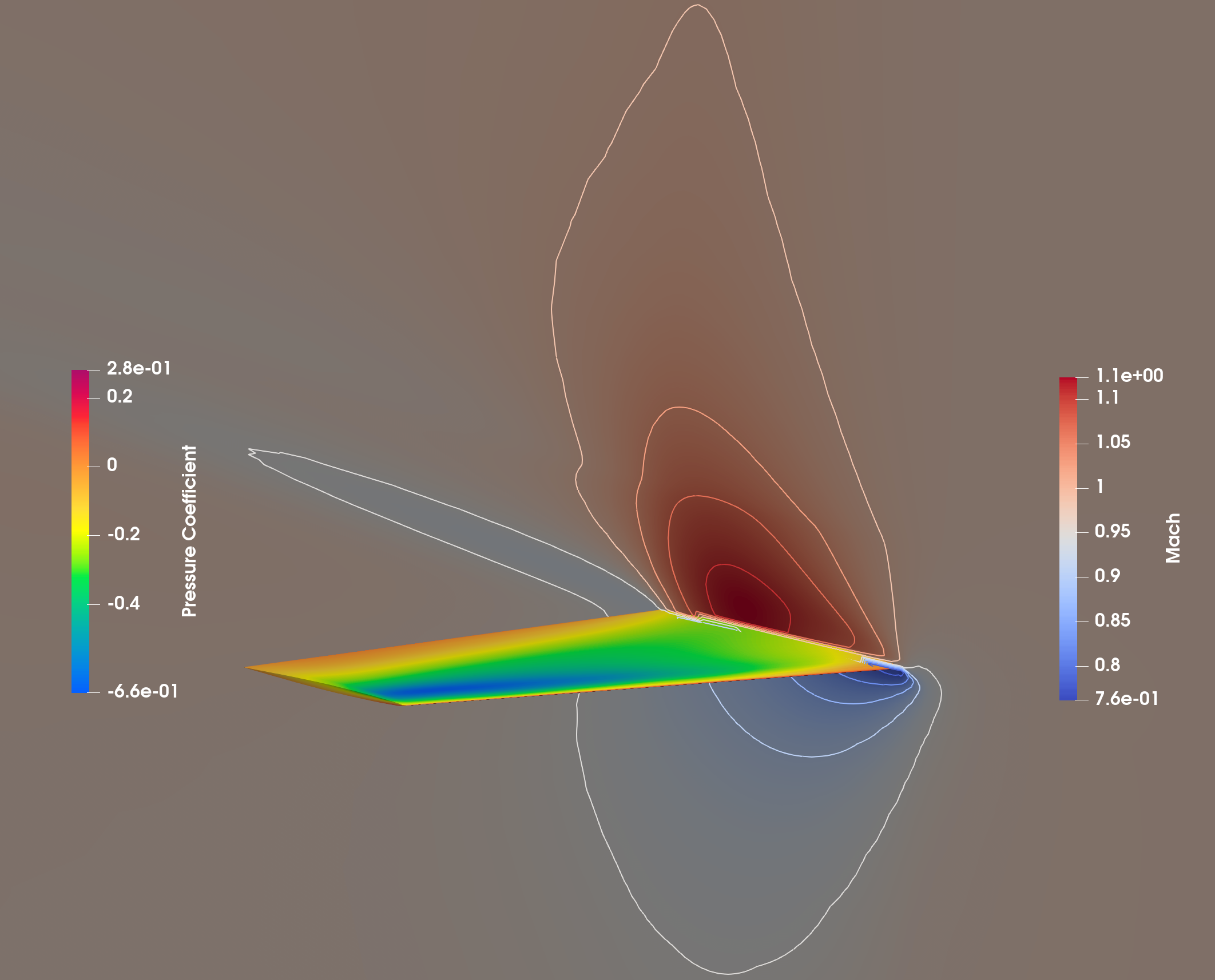}
    \caption{Initial condition for the FSI simulation: pressure coefficient on the wing and Mach contour profile in the fluid. The initial condition is also the converged steady state result with a fixed AGARD wing.}
    \label{fig:Wing-steady}
\end{figure}

\begin{figure}[ht]
\centering
    \includegraphics[width=0.24\textwidth]{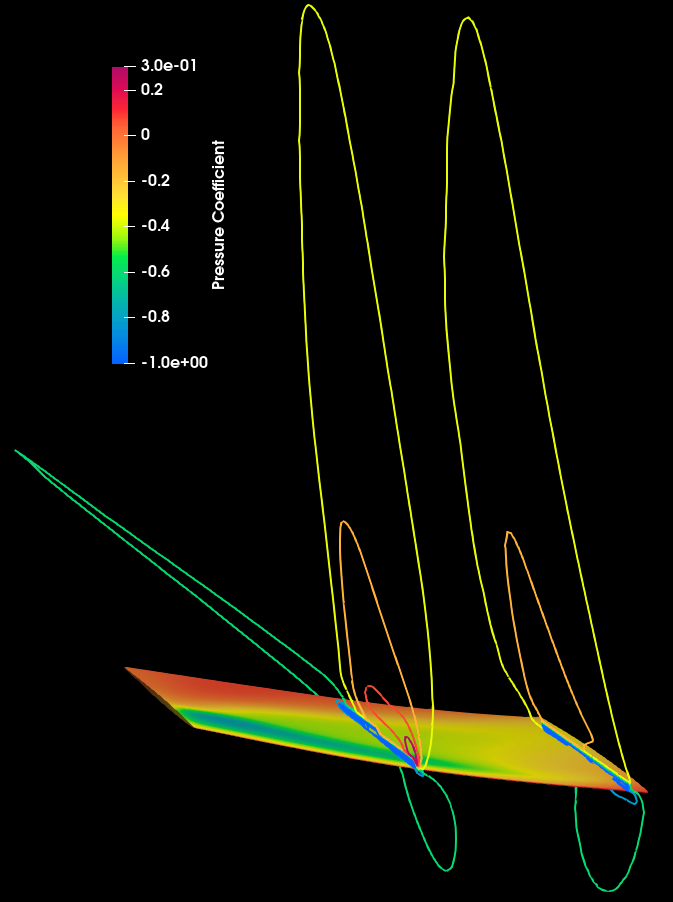}
    \includegraphics[width=0.24\textwidth]{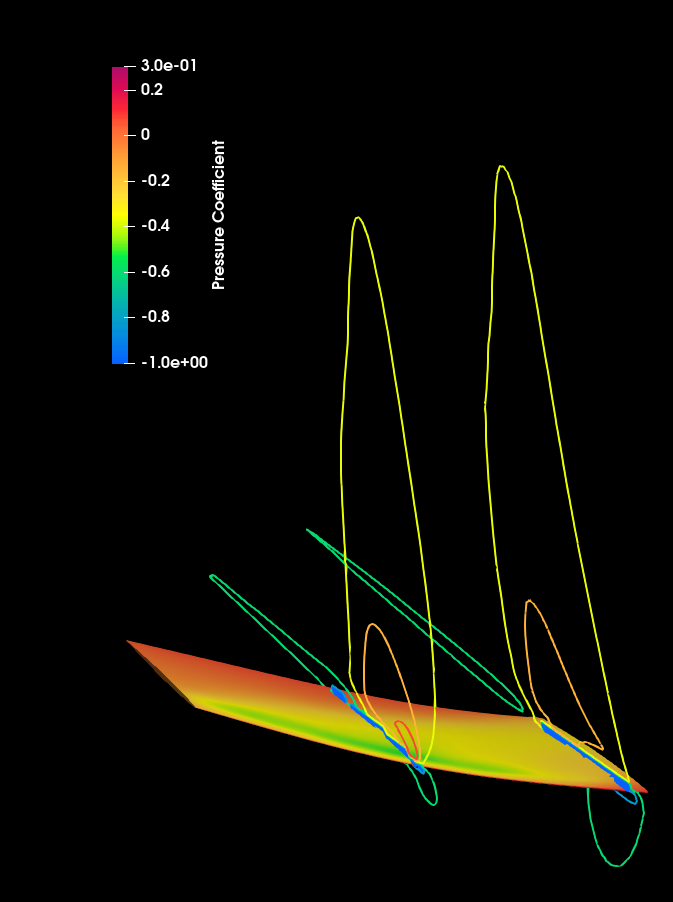}
    \includegraphics[width=0.24\textwidth]{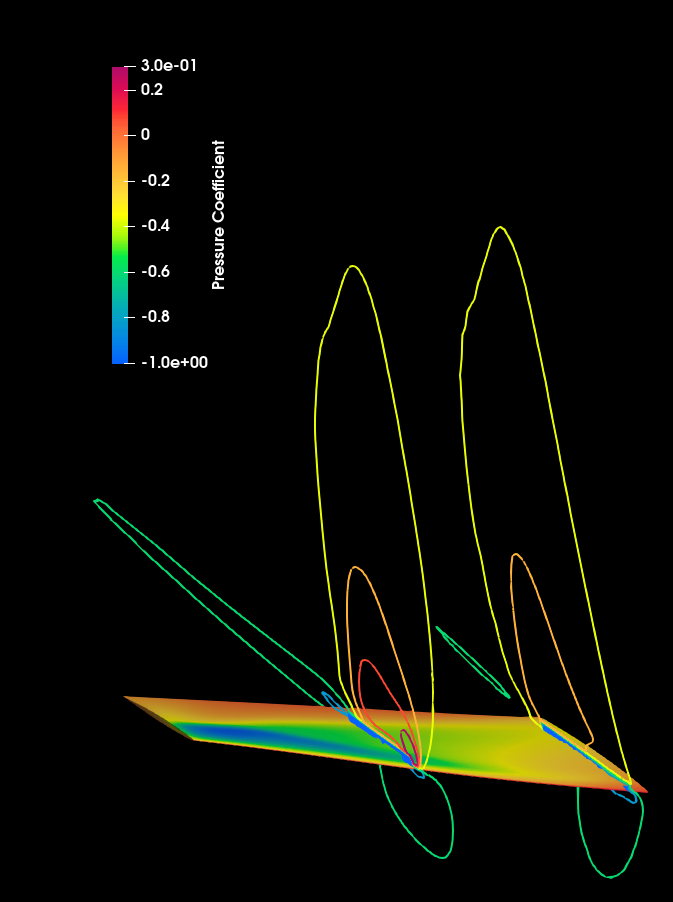}
    \includegraphics[width=0.24\textwidth]{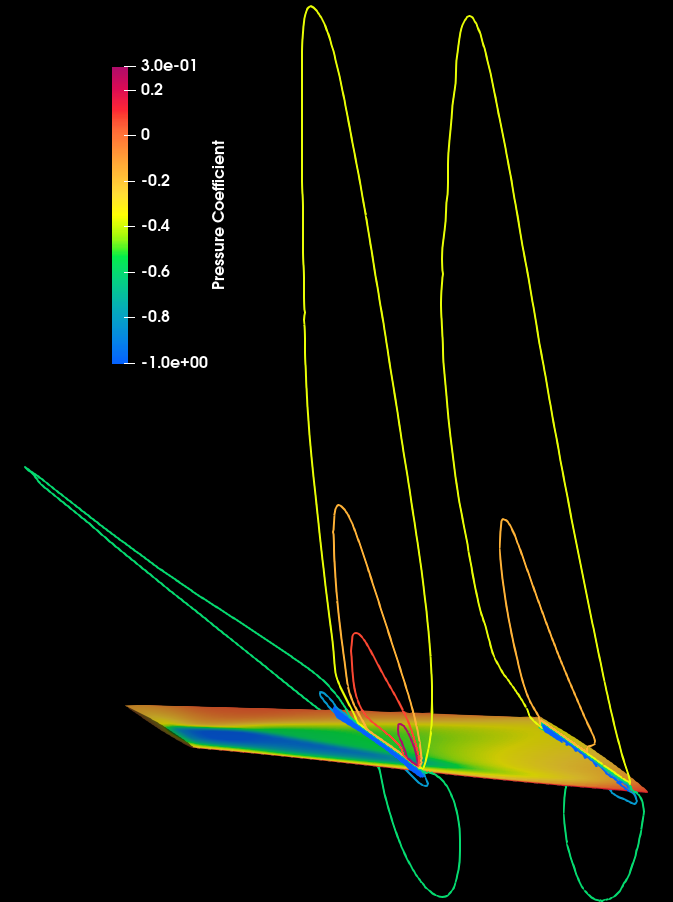}
    \caption{AGARD wing transonic buffeting: Snapshots of the fluid-structure solution at t = 0.025~s, t = 0.05~s, t = 0.075~s, t = 0.10~s (left to right) with pressure coefficient on the wing and Mach contour in the fluid.}
    \label{fig:Wing-unsteady}
\end{figure}

For the forward problem, we apply FIVER (Section~\ref{subsec:computational_framework}) to solve for the fluid dynamics and structural response. Specifically, an implicit backward Euler time integrator is used for the fluid and an implicit 2nd order mid-point rule time integrator is applied to for the structure with a constant time step $\Delta t = 5\times10^{-5}s$. 
The simulation starts from a steady flow field depicted in~\cref{fig:Wing-steady}, in which a shock wave forms on the wing.  When the coupled simulation starts, the wing starts buffeting, interacts with the shock wave, and undergoes large deformations~(See~\cref{fig:Wing-unsteady}). 
The coupled system is integrated till $T = 0.1s$. And each forward run takes about 36 core hours.

\begin{figure}[ht]
\centering
    \includegraphics[width=0.95\textwidth]{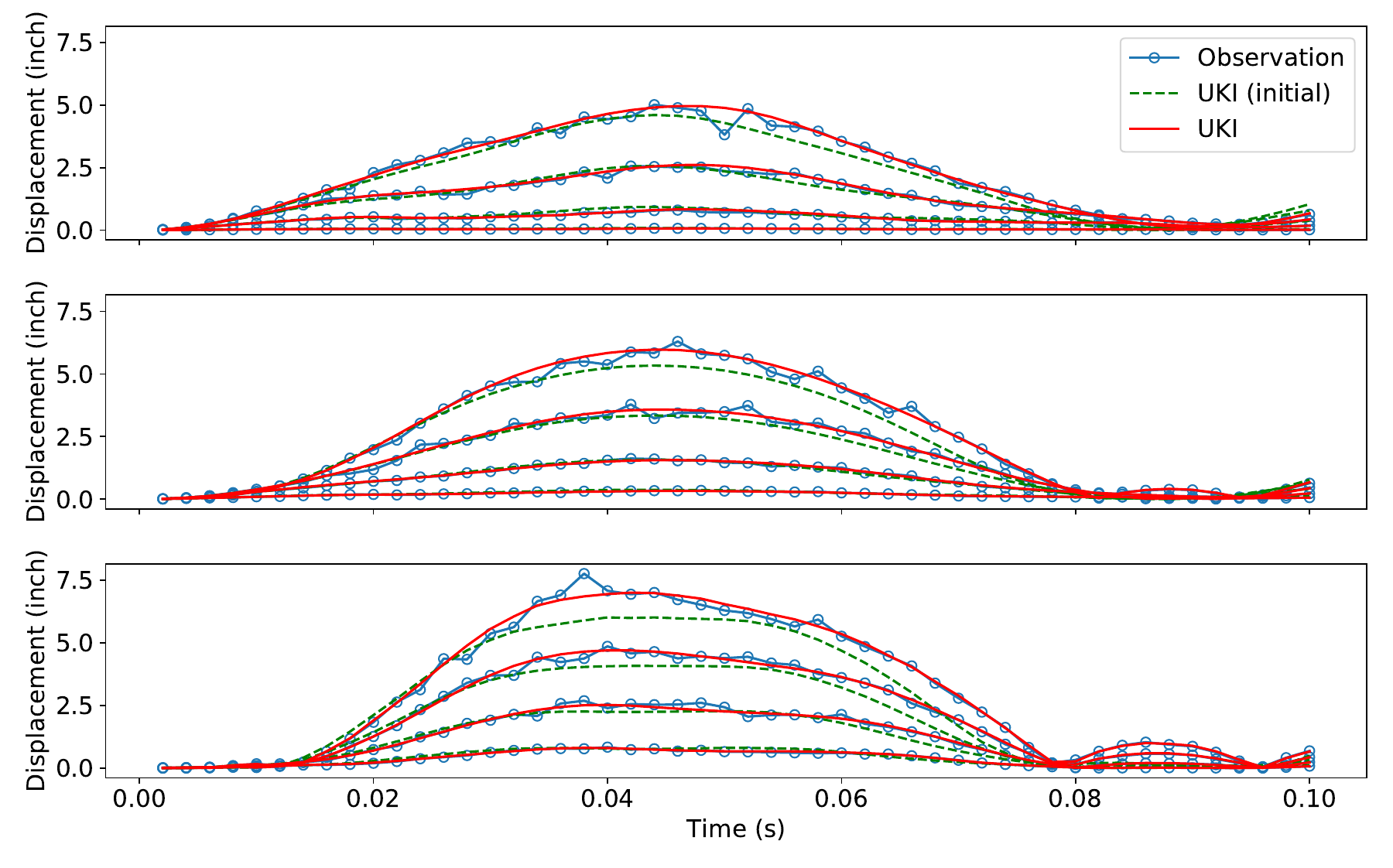}
    \caption{Blue dot lines: displacement observation data on the leading edge measurement locations (top), 
    on the half-chord measurement locations (middle), 
    on the trailing edge measurement locations (bottom) with the damaged AGARD wing. 
    Dashed green lines: predictions with the initial guess (no damage).
    Red lines: predictions with the damage field obtained by the UKI at the 15th iteration.
    }
    \label{fig:Wing-Disp}
\end{figure}

For the inverse problem, the damage field $\omega(y, \theta)$ is inferred from the displacement measurements of the wing buffeting motion. 
With the reference damage field, we construct the observation $y_{obs}$ at 12 locations~(See \cref{fig:Wing-Damage-Obs}) on the leading edge, half chord position, and trailing edge with
$5\%$ synthetic Gaussian random noises:
\begin{equation*}
    y_{obs} = y_{ref} + 5\% y_{ref} \odot \N(0, \I),
\end{equation*}
where $\odot$ denotes element-wise multiplication.
The observation data are collected every $0.002$~s, giving $N_y = 600$ observation data in total (See \cref{fig:Wing-Disp}).

\begin{figure}[ht]
\centering
    \includegraphics[width=0.6\textwidth]{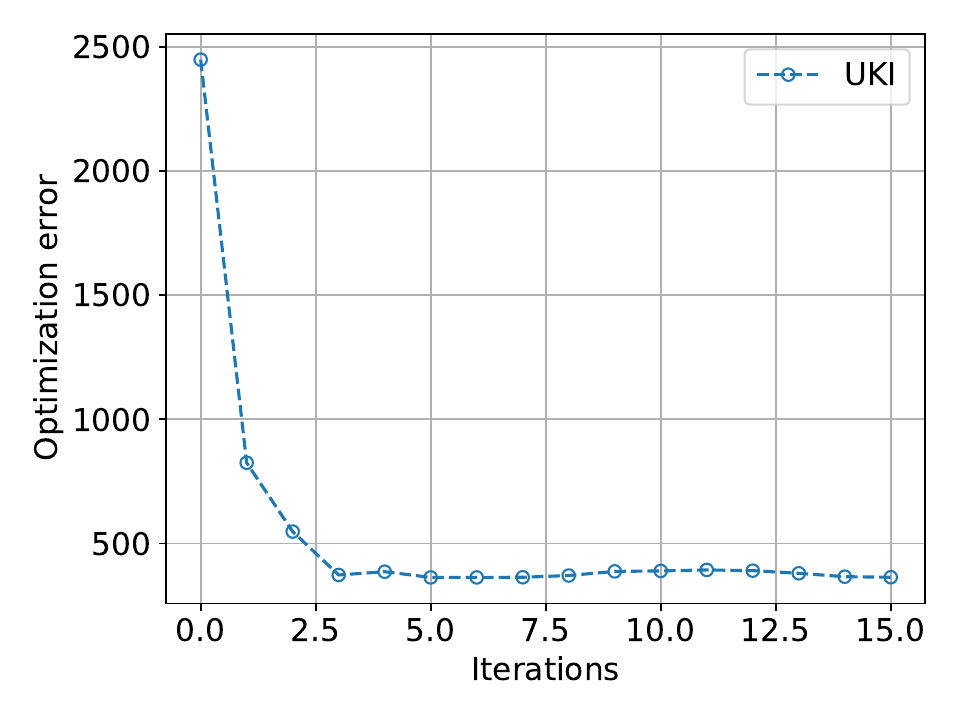}
    \caption{The optimization error $\frac{1}{2}\lVert\Sigma_{\eta}^{-\frac{1}{2}}(y - \G(\theta)) \rVert^2$ at each UKI iteration for the wing damage detection problem.}
    \label{fig:Wing-Error}
\end{figure}

\begin{figure}[ht]
\centering
    \includegraphics[width=0.6\textwidth]{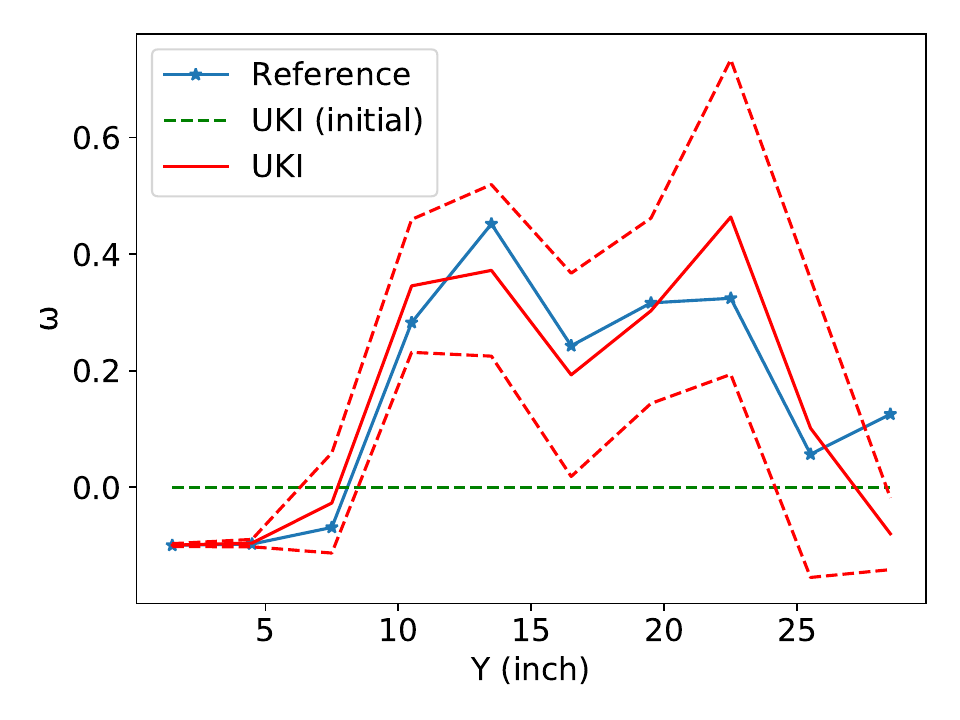}
    \caption{The estimated damage field $\omega(y, \theta_n)$~(red solid line) and the associated 2-$\sigma$ confidence intervals~(red dashed lines) at the $15$-th iteration; The reference damage field $\omega_ref(y)$~(blue star line); The initial guess of the damage field at the $0$-th iteration (green dashed line) for the wing damage detection problem.}
    \label{fig:Wing-UQ}
\end{figure}

The UKI is applied to retrieve the first $N_{\theta} = 5$ KL modes; in essence the model form error exists. The solution process is initialized with $\mean_0 = 0$ (no damage) and $\Cov_0 =  \I$, and hence $\theta_0 \sim \N(0, \I)$.
The observation error is set to follow $\eta \sim \N(0, 0.1^2\I)$; in essence we assume imperfect knowledge of the noise model.
For $N_{\theta} = 5$, each iteration includes $11$ parallel forward runs.
The optimization errors $\Phi(\theta)$ in \cref{eq:LSQ} at each iteration are plotted in~\cref{fig:Wing-Error}, showing fast convergence.
The estimated damage field $\omega(y, \theta)$ and the associated 2-$\sigma$ confidence intervals at the $15$-th iteration are depicted in \cref{fig:Wing-UQ}. The truth damage field $\omega_{ref}$ falls within the confidence interval with high probability. The predicted displacement fields at these measurement locations and the displacement predicted with the initial guess (no-damage case) are depicted in~\cref{fig:Wing-Disp}. It is worth mentioning that all displacement curves are very close to the observation data, which indicates the wing displacement is not very sensitive to the damage. However, UKI still delivers a good estimate of the damage field, which demonstrates the effectiveness of the Bayesian calibration procedure for real-world applications.

\section{Conclusion}
\label{sec:Conclusion}
This paper presents a general Bayesian calibration framework for FSI problems based on unscented Kalman inversion. It is attractive for at least four reasons: i) It is non-intrusive and derivative-free; ii) It is robust for chaotic inverse problems with noisy observations; iii) It provides uncertainty information; iv) It is embarrassingly parallel. It is well-adapted to parameter/field estimation problems for any large complex computational models given as a black box. There are numerous directions for future improvements:

\begin{itemize}
    \item At each iteration, UKI requires solving the forward problem  $2N_\theta+1$ times. Although these forward solver evaluations are embarrassingly parallel, the computational cost might be intractable when the number of parameters $N_{\theta}$ is large or the forward evaluation is expensive. The use of reduced-order models, including but not limit to neural network surrogate models~\cite{raissi2019physics,huang2020learning,xu2021learning,jimenez2019fatigue,kovachki2021multiscale}, projection-based reduced order models~\cite{dowell2001modeling,lieu2006reduced, balajewicz2014reduction,xiao2016non,taira2017modal}, and Kernel-based surrogate models~\cite{quinonero2005unifying,hofmann2008kernel,hamzi2021data}, is worth exploring in the future.
    \item For the damage detection problem in Subsection~\ref{ssec:damage}, the sensors are uniformly located on the aircraft wing. Optimal sensor placement~\cite{worden2001optimal,meo2005optimal,manohar2018data}, which potentially makes data collection more efficient and makes the algorithm converge faster, is worth further investigation. And more realistic damage and fracture models~\cite{moes1999finite,zhan2017discontinuous,zhan2018new} will be considered in the future.
    \item For Bayesian calibration problems with time-series data studied in the present work, other Kalman methodology~(e.g. Kalman smoother~\cite{evensen2000ensemble,bocquet2014iterative,spratt2021characterizing}) will be applied and studied in the future.   
\end{itemize}

\section{acknowledgement*}
The authors gratefully acknowledge the support of National Institutes of Health under Award P01-DK043881 and the generosity of Eric and Wendy Schmidt by recommendation of the Schmidt Futures program. 
The authors thank Dr. Kevin G. Wang and Dr. Fangbao Tian for their advice on the manuscript.

\section{data availability statement*}
All computer code used in this paper is open source. 
Datasets, including mesh files and results, are available at \url{https://github.com/Zhengyu-Huang/InverseProblems.jl}.

\appendix

\section{Proof of Theorem~\ref{TH:LINEAR}}
\label{sec:linear-proof}
Thanks to the linearity, the Eqs.~\eqref{eq:KF_joint2} are reduced to 
\begin{align*}
    \py_{n+1} = G\mean_n, \quad \pCov_{n+1}^{\theta p} = \pCov_{n+1} G^T,\quad  \textrm{and} \quad \pCov_{n+1}^{pp} = G  \pCov_{n+1} G^T + \Sigma_{\nu}.
\end{align*}
The update equations~\eqref{eq:KF_analysis} become
\begin{equation}
\label{eq:Lin_KF_analysis}
    \begin{split}
        \mean_{n+1} &= \mean_{n} + \pCov_{n+1} G^T (G  \pCov_{n+1} G^T + \Sigma_{\nu})^{-1} (y - G\mean_{n}), \\
         \Cov_{n+1}&= \pCov_{n+1} - \pCov_{n+1} G^T(G  \pCov_{n+1} G^T + \Sigma_{\nu})^{-1} G \pCov_{n+1},
    \end{split}
\end{equation}
with $\pCov_{n+1} =  \Cov_{n} + \Sigma_{\omega}$.

With the hyperparameters defined in~\cref{eq:hyperparameters}, the update equation of $\{\Cov_n\}$ in \cref{eq:Lin_KF_analysis} can be rewritten as 
\begin{equation}
\label{eq:Lin_KF_Cinv}
\begin{split}
    &\Cov_{n+1}^{-1} = G^T\Sigma_{\nu}^{-1}G + (\Cov_n + \Sigma_{\omega})^{-1} = \frac{1}{2}G^T\Sigma_{\eta}^{-1}G + (2\Cov_n)^{-1}.\\
\end{split}
\end{equation}
We have a closed formula for $\Cov^{-1}_{n}$: 
\begin{equation}
\label{eq:beta_1}
    \Cov_n^{-1} = 
    \Big[1 - \frac{1}{2^n}\Big] G^T\Sigma_{\eta}^{-1}G + \frac{1}{2^n} \Cov_0^{-1}.
\end{equation}
This leads to the exponential convergence $\displaystyle \lim_{n\to \infty} \Cov_n^{-1} = G^T\Sigma_{\eta}^{-1}G$.

The convergence proof of $\mean_n$ basically follows the work of Huang {\it{et al.}}~\cite{UKI1}. \Cref{eq:Lin_KF_Cinv,eq:beta_1} lead to 
\begin{equation}
\label{eq:C_bound}
\begin{split}
\frac{1}{2}G^T\Sigma_{\eta}^{-1}G = G^T\Sigma_{\nu}^{-1}G \preceq\Cov_{n+1}^{-1} \preceq G^T\Sigma_{\nu}^{-1}G + \Sigma_{+}
    \quad \textrm{where} \quad 
    \Sigma_{+} = 
     G^T\Sigma_{\nu}^{-1}G + \Cov_0^{-1}. 
\end{split}
\end{equation}
The update equation of $\mean_n$ in \cref{eq:Lin_KF_analysis} can be rewritten as 
\begin{equation}
\label{eq:contracting-1}
\mean_{n+1} = \mean_{n} + \Cov_{n+1}G^T\Sigma_{\nu}^{-1}(y - G\mean_{n}).
\end{equation}
We have the assumption that $G$ has full column rank, and therefore, $B := G^T\Sigma_{\nu}^{-1}G$ is positive definite. From this, it follows that $\I -\Cov_{n+1} B$ has the same spectrum as $\I - B^{\frac{1}{2}}\Cov_{n+1}B^{\frac{1}{2}}$. Using the bounds on $\Cov_{n+1}$ appearing in~\cref{eq:C_bound}, the spectral radius of the update matrix in \cref{eq:contracting-1} satisfies
\begin{equation}
\label{eq:cov-mean-par}
\begin{split}
\rho(\I - \Cov_{n+1}G^T\Sigma_{\nu}^{-1}G) 
&= \rho(\I - \Cov_{n+1}B) \\
&= \rho\Big(\I - \sqrt{B}\Cov_{n+1}\sqrt{B}\Big) \\
&\leq
1 - \rho\Big(\sqrt{B}\big(B + \Sigma_{+}\big)^{-1}\sqrt{B}\Big) \\
&= 1 -\epsilon_0,
\end{split}
\end{equation}
where $\epsilon_0 \in (0,1)$. Hence, we have that $\{\mean_{n}\}$ converges exponentially to $\mean_{\infty}$, which satisfies $\displaystyle G^T\Sigma_{\nu}^{-1}(y - G\mean_{\infty}) = 0$. And it is a minimizer of $\Phi$.

As for the uncertainty estimation, we have 
\begin{equation*}
    G^T \Sigma_{\eta}^{-1}(y - G\mean_{\infty}) = 0 \qquad
    \Cov^{-1}_{\infty} = G^T\Sigma_{\eta}^{-1}G.
\end{equation*}
and 
\begin{equation*}
    \eta_r = y - G\theta_{ref} \sim \N(0, \Sigma_{\eta}).
\end{equation*}
When $G$ is of full column rank, $G^T\Sigma_{\eta}^{-1}G$ is non-singular. We have
\begin{equation*}
\begin{split}
&\theta_{ref} - \mean_{\infty} =  - \Big(G^T\Sigma_{\eta}^{-1}G \Big)^{-1} G^T \Sigma_{\eta}^{-1}\eta_r ,\\
&\theta_{ref} - \mean_{\infty} \sim \N\Big(0, \big(G^T\Sigma_{\eta}^{-1}G\big)^{-1}\Big) = \N(0, \Cov_{\infty}) .
\end{split}
\end{equation*}
Therefore, each component of $\theta_{ref} - \mean_{\infty}$ obeys Gaussian distribution, 
\begin{equation*}
    {\theta_{ref}}_{(i)} - {\mean_{\infty}}_{(i)} \sim \N(0, {\Cov_{\infty}}_{(i,i)}).
\end{equation*}
The empirical rule of Gaussian implies \Cref{eq:err_bound_linear}.

\bibliography{references}{}
\end{document}